\documentclass[12pt,twoside]{article} 
\usepackage[latin1]{inputenc}
\usepackage{graphicx}
\usepackage{enumitem}
\usepackage{latexsym}
\usepackage{amssymb}
\usepackage{epsfig}
\usepackage{xcolor}
\usepackage{amsmath}
\usepackage{float}
\usepackage{caption}
\usepackage{xr}

\makeatletter

\@addtoreset{equation}{section}
\makeatother

\allowdisplaybreaks 

\def\eps{\varepsilon}

\def\e{\varepsilon}

\def\beq{\begin{eqnarray*}}
\def\eeq{\end{eqnarray*}}

\def\Lt{\Lambda_\vartheta}
\def\tin{x_{i,n}}
\def\t{\vartheta}

\newcommand{\nto}{\xrightarrow[n\to\infty]{}}

\newtheorem{theo}{Theorem}[section]
\newtheorem{lemma}[theo]{Lemma}
\newtheorem{cor}[theo]{Corollary}
\newtheorem{rem}[theo]{Remark}

\newtheorem{example}[theo]{Example}

\begin{document}

\title{\bf Semi-parametric transformation boundary regression models\footnote{Financial support by the DFG (Research Unit FOR 1735 {\it Structural Inference in
Statistics: Adaptation and Efficiency}) is gratefully acknowledged.}}

\author{{\sc Natalie Neumeyer, Leonie Selk} and {\sc Charles Tillier}\\ Department of Mathematics, University of Hamburg}

\maketitle

\begin{abstract}
In the context of nonparametric regression models with one-sided errors, we consider parametric transformations of the response variable in order to obtain independence between the errors  and the covariates. 
In view of estimating the tranformation parameter, we use a minimum distance approach and show the uniform consistency of the estimator under mild conditions. The boundary curve, i.e.\ the regression function, is estimated applying a smoothed version of a local constant approximation for which we also prove the uniform consistency. We deal with both cases of random covariates and deterministic (fixed) design points. To highlight the applicability of the procedures and to demonstrate their performance, the small sample behavior is investigated in a simulation study using the so-called Yeo-Johnson transformations. 
\end{abstract}






\section{Introduction} 

Before fitting a regression model it is very common in applications to transform the response variable. The aim of the transformation is to gain efficiency in the statistical inference, for instance, by reducing skewness or inducing  a specific structure of the model, e.g.\  linearity of the regression function or homoscedasticity.
 In practice often a parametric class of transformations
is considered from which an `optimal' one should be selected data-dependently (with a specific purpose in mind).  A classical example is the class of Box-Cox power transformations introduced for linear models by Box and Cox (1964).
There is a vast literature on parametric transformation models in the context of mean regression and we refer to the monograph by Carroll and Ruppert (1988).
Powell (1991) introduced Box-Cox transformations in the context of linear quantile regression; see also
Mu and He (2007) who considered transformations to obtain a linear quantile regression function.
Horowitz (2009) reviewed estimation in transformation models with parametric regression in the cases where either the transformation or the error distribution or both are modeled nonparametrically.
Linton, Sperlich and Van Keilegom (2008) suggested parametric estimators for transformations, while the error distribution is estimated nonparametrically and the regression function is additive.
In this paper, the aim of the transformation is to induce independence between the covariables and the errors. Linton et al.\ (2008) considered profile likelihood and minimum distance estimation for the transformation parameter. 
The results for the profile likelihood estimator were generalized for nonparametric regression models by Colling and Van Keilegom (2016).\newline
\hspace*{0.6cm} All literature cited above is about mean or quantile regression. In contrast in the paper at hand we consider boundary regression models. Such nonparametric regression models with one-sided errors have been considered, among others, by Hall and Van Keilegom (2009), Meister and Rei{\ss} (2013), Jirak, Meister and Rei{\ss} (2014) and Drees, Neumeyer and Selk (2018). Relatedly, estimation of support boundaries have been considered, for instance, by Hõrdle, Park and Tsybakov (1995), Hall, Park and Stern (1998), Girard and Jacob (2008) and Daouia, Noh and Park (2016). Such models naturally appear when analyzing auctions or records or production frontiers. 
Unlike conditional mean models, regression models with one-sided errors (as well as quantile regression models)
have the attractive feature of equivariance under monotone transformations. Thus in such a model with monotone transformation of the response one can recover the original functional dependence in an easy manner.
Similar to Linton et al.\ (2008) the aim of our transformation is to induce a model where the error distribution does not depend on the covariates. 
Independence of errors and covariates is a very typical assumption in regression models. For boundary models this assumption is met, e.g., by M³ller and Wefelmeyer (2010), Meister and Rei{\ss} (2013), and Drees et al.\ (2018). 
A transformation inducing (approximate) independence between the covariable and the  error would allow for a global bandwidth selection in the adaptive regression estimator suggested by Jirak et al.\ (2014). 
Wilson (2003) pointed out that in production frontier models, independence assumptions are needed for validity of  bootstrap procedures for nonparametric frontier models (see
Simar and Wilson, 1998) and suggested some tests for independence of errors and covariates (see also Drees et al., 2018).

While Linton et al.\ (2008) found advantages of the profile likelihood approach over minimum distance estimation  of the transformation parameter in corresponding mean regression transformation models, this is at the cost of strong regularity conditions, among others a bounded error density with bounded derivative. In the context of boundary models with  error distribution which is regularly varying at zero and irregular, one needs to avoid assumptions on bounded densities. Thus we investigate a minimum distance approach to estimate the transformation parameter and give mild model assumptions under which the estimator is consistent.\newline 
\hspace*{0.6cm} We consider the cases of random covariates and deterministic (fixed) design points, which are both meaningful. The equidistant fixed design - as well as its natural generalization to deterministic covariates - is often used in real-life applications when time is involved in the data set. This is the case for instance in Jirak et al.\ (2014) where the authors studied the monthly sunspot observations and the annual best running times of 1500 meters. Besides, deterministic design is met accross a number of papers in regression models, see for instance Brown and Low (1996), Meister and Rei\ss \ (2013) and the references within. The case of random covariates is obviously the most relevant and appears in essence in many applications in boundary models, among other, in insurance and financial risk modelling when analyzing the optimality of portfolios (see Markowitz (1952) for the seminal contribution).

The remaining part of the manuscript is organized as follows. In section 2 the model is explained, while in section 3 the estimation procedure is described. In section 4 we show consistency of the  transformation parameter estimator. In section 5 we present simulation results. The proofs for the random covariate case are given in the appendix, while supplementary material contains proofs for the fixed design case and some additional figures and simulation results.

\section{Model}

\subsection{The random design case}

Consider independent and identically distributed observations $(X_i,Y_i)$, $i=1,\ldots,n$, with the same distribution as $(X,Y)$, where $Y$ is univariate and  $X$ is distributed on $[0,1]$.
Further consider a family  $\mathcal{L}=\{\Lt|\vartheta\in\Theta\}$ of strictly  increasing and continuous transformations.
Throughout the paper we assume existence of a transformation $\Lambda_{\vartheta_0}$ in the class $\mathcal{L}$ such that in the corresponding boundary regression model
\begin{equation}\label{model0-random}
\Lambda_{\vartheta_0}(Y)=h_{\vartheta_0}(X)+\e
\end{equation}
 the errors and the covariates are stochastically independent. Note that for notational simplicity we  set $\Lambda_0=\Lambda_{\vartheta_0}$ and $h_0=h_{\vartheta_0}$. Further denote by $F_0$ the cumulative distribution function (cdf) of the independent and identically distributed (iid) $\e_i=\Lambda_0(Y_i)-h_0(X_i)$, $i=1,\dots,n$. Then we assume that $F_0(0)=1$ and $F_0(-\Delta)<1$ for all $\Delta>0$. This identifies the function $h_0$ as the upper boundary curve of the observations since  
\begin{eqnarray*}
\mathbb{P}(\Lambda_0(Y_{i}) \leq h_0(X_i)\mid X_i=x)&=&1  \mbox{ for all } x\in[0,1]\\
\mathbb{P}(\Lambda_0(Y_{i}) -h_0(X_i)\leq -\Delta\mid X_i=x)&<&1 \mbox{ for all }x\in [0,1], \Delta>0. 
\end{eqnarray*} 
The aim is to estimate $\vartheta_0$ from the observations. 
\begin{rem}
Note that even if the model does not hold exactly (i.e.\ there does not exist any $\vartheta_0\in \Theta$ that leads to exact independence of the errors and covariates) the transformation can be useful in applications because it will reduce the dependence. 
\end{rem}


For each $\vartheta\in\Theta$ one can consider the transformed responses $\Lambda_\vartheta (Y_{i})$. Note that those form a boundary regression model with boundary curve $ h_{\vartheta}=\Lt\circ \Lambda_0^{-1}\circ h_0$, because
$$\mathbb{P}(\Lambda_\vartheta(Y_{i})\leq h_\vartheta(X_i)\mid X_i=x)=\mathbb{P}(\Lambda_0(Y_{i})\leq h_0(X_i)\mid X_i=x)=1$$ 
and for each $\delta>0$,
$$\mathbb{P}\left(\Lambda_\vartheta(Y_{i})-h_\vartheta(X_i)\leq -\delta\mid X_i=x\right)=\mathbb{P}\left(\Lambda_0(Y_{i})\leq \Lambda_0(\Lambda_\vartheta^{-1}\left(h_\vartheta(x)-\delta)\right)\mid X_i=x\right)<1$$
since $\Delta=h_0(x)-\Lambda_0(\Lambda_\vartheta^{-1}(\Lt(\Lambda_0^{-1}(h_0(x)))  -\delta))>0$ since each $\Lambda_\vartheta$ is strictly increasing.
The conditional distribution of $\Lt(Y_{i})$ for some general $\vartheta\in\Theta$  reads as
\beq
\mathbb{P}(\Lt(Y_{i})\leq y\mid X_i=x)&=&
 \mathbb{P}\left(\Lambda_0(Y_{i})\leq\Lambda_0(\Lt^{-1}(y))\mid X_i=x\right)
\;=\;F_0\left(\Lambda_0(\Lt^{-1}(y))-h_0(x)\right).
\eeq

\begin{rem}\label{identifiability}
It is important to give conditions under which the unknown components $\Lambda_0=\Lambda_{\vartheta_0}$, $h_0=h_{\vartheta_0}$ and $F_0$ in model (\ref{model0-random}) are identifiable. To this end we impose the following assumptions. 
\begin{itemize}
\item Assume that $Y$ has a continuous distribution and w.l.o.g.\ assume that 0 is in the data range (otherwise shift the data). 
\item Assume that $X$ is continuously distributed with support $[0,1]$.
\item Assume $\Lambda_\vartheta(0)=0$ for all $\vartheta\in\Theta$, and $\Lambda_\vartheta$ is strictly increasing and continuous for each $\vartheta\in\Theta$.
\item Assume that if for some $\vartheta_0,\vartheta_1\in\Theta$ one has
$$(\Lambda_{\vartheta_1}\circ\Lambda_{\vartheta_0}^{-1})(a-b)=(\Lambda_{\vartheta_1}\circ\Lambda_{\vartheta_0}^{-1})(a)-(\Lambda_{\vartheta_1}\circ\Lambda_{\vartheta_0}^{-1})(b)$$
for all $a,b\in J$, where $J$ is an interval of positive length, then it follows that $\vartheta_0=\vartheta_1$. 
\item Assume that $F_0$ (the cdf of $\varepsilon=\Lambda_{\vartheta_0}(Y)-h_{\vartheta_0}(X)$) is strictly increasing.
\item Assume that $h_{\vartheta_0}$ is not constant and is continuous. 
\end{itemize}
Now assume that the model 
$$\Lambda_\vartheta(Y)=h_\vartheta(X)+\varepsilon(\vartheta)\mbox{ with } X \mbox{ independent from }\varepsilon(\vartheta)$$
holds for $\vartheta=\vartheta_0$ (with our notations $\Lambda_{\vartheta_0}=\Lambda_0$, $h_{\vartheta_0}=h_0$, $\varepsilon(\vartheta_0)=\varepsilon$) and for $\vartheta=\vartheta_1$.

Note that from the assumption it follows that $\mathbb{P}(\varepsilon(\vartheta)\leq 0)=1$, $\mathbb{P}(\varepsilon(\vartheta)\leq -\Delta)<1$ for each $\Delta>0$, such that $h_\vartheta$ is the upper boundary curve in the model (for $\vartheta\in \{\vartheta_0,\vartheta_1\}$).

We show in section \ref{ident-proof} of the appendix that it follows that $\vartheta_0=\vartheta_1$. Thus the transformation is identifiable. Further $h_{\vartheta_0}(x)$ is then the right endpoint of the conditional distribution of $\Lambda_{\vartheta_0}(Y)$, given $X=x$, and $F_0$ is identified as cdf of $\Lambda_{\vartheta_0}(Y)-h_{\vartheta_0}(X)$. 

If the function class $\mathcal{L}$ contains the identity, then the assumptions rule out that it contains transformations which are linear on some interval with positive length. On the other hand it is clear that linear transformations can  never reduce the dependence between the covariate and  the error distribution. 
\end{rem}

\begin{example}\label{yeo-johnson} In this example we give two classes of transformations that fulfill the identifiability assumptions. 

 Yeo and Johnson (2000) generalized the Box-Cox transformations by suggesting
\begin{equation*}
\Lambda_\vartheta(y) = \begin{cases} \frac{(y+1)^\vartheta-1}{\vartheta}, &\text{if } y \ge 0,\vartheta\neq 0 \\
 \log(y+1), &\text{if } y\ge 0,\vartheta=0\\
 -\frac{(-y+1)^{2-\vartheta}-1}{2-\vartheta}, &\text{if } y<0,\vartheta\neq 2 \\
 -\log(-y+1), &\text{if } y<0,\vartheta=2,
 \end{cases}
 \end{equation*}
which are typically considered for $\t\in\Theta=[0,2]$ because then they are bijective maps $\Lambda_\t:\mathbb{R}\to \mathbb{R}$. 
Note that  $\Lambda_\vartheta(0)=0$ for all $\vartheta\in\Theta$.

The class of sinh-arcsinh transformations, see Jones and Pewsey (2009), do shift the location, but they can be modified to fulfill $\Lambda_\vartheta(0)=0$ for all $\vartheta\in\Theta$, e.g. consider
$$\Lambda_{(\vartheta_1,\vartheta_2)}(y)=\sinh (\vartheta_1 \sinh^{-1}(y)-\vartheta_2)-\sinh(-\vartheta_2).$$
Here $\vartheta_1>0$ is the tailweight parameter and $\vartheta_2 \in \mathbb{R}$ the skewness parameter. These transformations define also bijective maps  $\Lambda_{(\vartheta_1,\vartheta_2)}:\mathbb{R}\to \mathbb{R}$. 

\end{example}


\subsection{The fixed design case}
In the fixed design case we consider a triangular array of independent observations $Y_{i,n}$, $i=1,\ldots,n$, and deterministic design points $0<x_{1,n}<\dots <x_{n,n}<1$. 
Once again we assume existence of a transformation $\Lambda_0=\Lambda_{\vartheta_0}$ in the class $\mathcal{L}$ such that setting $h_0=h_{\vartheta_0}$ in the corresponding regression model
\begin{equation}\label{model0}
\Lambda_0(Y_{i,n})=h_0(\tin)+\e_{i,n}
\end{equation}
the cdf of the errors does not depend on the design points, i.e.\ $\e_{i,n}\sim F_0$ $\forall i,n$. Note that, as in the random design case, we assume $F_0(0)=1$ and $F_0(-\Delta)<1$ for all $\Delta>0$ leading again to $ h_{\vartheta}=\Lt\circ \Lambda_0^{-1}\circ h_0$. 

\begin{rem}\label{identifiability-fixed}
Identifiability can be shown under the same conditions as in Remark \ref{identifiability}
as long as $ \bar{\Delta}_n:=\max_{1 \leq i \leq n+1} \left(x_{i,n}-x_{i-1,n}\right)\to 0 $; see section \ref{ident-proof-fixed} of the supplement. 
\end{rem}

\begin{example}\label{example1}Figures \ref{graphic1} and \ref{graphic2} show realizations of the original data and the transformed data \eqref{model0} using a Yeo and Johnson transformation; see Example \ref{yeo-johnson}. For each figure, in the upper left panel the original data $(x_{i,n},Y_{i,n})$, $i=1,\dots,n=100$, are depicted with their boundary curve, while the upper right panel shows the corresponding non-iid errors. The lower left panel shows the transformed data with the curve $h_0$,  while the lower right panel shows the iid errors $\eps_{i,n}$, $i=1,\dots,n$.
\end{example}

\section{Estimating the transformation}

\subsection{The random design case}

If $\vartheta_0$ were known we could estimate the regression function (upper boundary curve) $h_0$ by a local constant approximation, i.e.
\begin{equation}\label{h_0}
\tilde h_0(x)=\max\{\Lambda_0(Y_{i})|i=1,\ldots,n\text{ with }|X_i-x|\leq b_n\},
\end{equation}
where $b_n\searrow 0$ is a sequence of bandwidths. For this estimator we will show uniform consistency 
under the following assumptions. 

\begin{enumerate}[label=(\textbf{A\arabic{*}})]
\item \label{F1} Model (\ref{model0-random}) holds with iid $\eps_{1},\dots,\eps_{n}\sim F_0$ and $F_0(0)=1$, $F_0(-\Delta)<1$ for all $\Delta>0$, and $\eps_{1},\dots,\eps_{n}$ are independent of $X_{1},\dots,X_{n}$.  
\item \label{X1} The covariates $X_1,\dots,X_n$ are iid with  cdf $F_X$ and density $f_X$ that is continuous and bounded away from zero on its support $[0,1]$. 
\item \label{H1} The regression function $h_0$ is continuous on $[0,1]$.
\item \label{B1} Let $(b_n)_{n\in\mathbb{N}}$ be a sequence of positive bandwidths that satisfies
 $\lim_{n\to\infty}b_n=0$ and $\lim_{n\to\infty}(\log n)/(nb_n)=0$.
\end{enumerate}

Note that we do not require any assumption on the error distribution. In particular, in the setup of regularly varying distributed errors, all the results hold for regular as well as irregular distributions. In what follows, let $\|\cdot\|_{\infty}$ denote the supremum norm and $I\{\cdot\}$ the indicator function. 

\begin{lemma}\label{consistency-random} Under model  (\ref{model0-random}) with assumptions \ref{F1}--\ref{B1} we have $\|\tilde h_0-h_0\|_{\infty}=o_P(1)$. 
\end{lemma}

The proof of the lemma is given in section \ref{proof-consistency-random} of the appendix. 
The result applies for a model without transformation. Thus, as a by-product, we show uniform consistency of a boundary curve estimator in models with random covariates (and non-equidistant fixed design, see Lemma \ref{consistency}), while in contrast Drees et al.\ (2018) assumed equidistant design and obtained rates of convergence under stronger assumptions on the error distribution $F_0$ and on the boundary curve $h_0$.

\noindent For general $\vartheta\in\Theta$ we define a simple boundary curve  estimator accordingly as
\[\tilde h_{\vartheta}(x)=\max\{\Lambda_{\vartheta}(Y_{i})|i=1,\ldots,n\text{ with }|X_i-x|\leq b_n\}\]
and it holds that $\tilde h_{\vartheta}=\Lt\circ \Lambda_0^{-1}\circ \tilde h_0$. Thus $\tilde h_\t$ consistently estimates $h_\t$. 
The local constant estimator can be improved by introducing slight smoothing. To this end, let $K$ be a density with compact support and $a_n$ some sequence of bandwidths that decreases to zero such that $na_n\to \infty$. 
Define
\begin{equation}\label{hatK}
\hat h_{\vartheta}(x)=\frac{\sum_{i=1}^n \tilde h_{\vartheta}(X_i)K\left(\frac{x-X_i}{a_n}\right)}{\sum_{i=1}^n K\left(\frac{x-X_i}{a_n}\right)},
\end{equation}
then $\hat h_{\vartheta}$ is also uniformly consistent for $h_\vartheta$; see Lemma \ref{lemma-smoothing}. 

\begin{example}For data as in Example \ref{example1}, Figures \ref{graphic3} and \ref{graphic4} in the online supplementary material demonstrate the smoothing of the estimator. We use the Epanechnikov-kernel $K(x)=0.75(1-x^2)I_{[-1,1]}(x)$ and bandwidths $b_n=0.5n^{-1/3}$, $a_n=0.5b_n$ with $n=100$. 
\end{example}

Based on this estimator we define 
the joint empirical distribution function of residuals and covariates as
$\hat F_{n,\vartheta}(y,s)=\frac 1n\sum_{i=1}^n I\{\Lt(Y_{i})-\hat h_\t(X_i)\leq y\}I\{X_i\leq s\}$.
For $\t=\t_0$, the covariate $X_i$ and the error $\Lt(Y_{i})- h_\t(X_i)$ are stochastically independent and thus, the joint empirical distribution function minus the product of the marginals, namely $\hat F_{n,\vartheta}(y,s)-\hat F_{n,\vartheta}(y,1)\hat F_{X,n}(s)$, estimates zero for $\t=\t_0$. Here $\hat F_{X,n}(\cdot)=\hat F_{n,\vartheta}(\infty,\cdot)$ denotes the empirical distribution function of $X_1,\dots,X_n$.  We will use this idea to estimate the transformation parameter $\t_0$. To this end, for any function $h:[0,1]\to\mathbb{R}$ define
\begin{eqnarray}\label{G_n-random}
 G_n(\t, h)(y,s)&=& \frac 1n\sum_{i=1}^n I\{\Lt(Y_{i})-h(X_i)\leq y\} \big( I\{X_i\leq s\}-\hat F_{X,n}(s)\big)
\end{eqnarray}
and note that $G_n(\t,\hat h_\t)(y,s)=\hat F_{n,\t}(y,s)-\hat F_{n,\t}(y,1)\hat F_{X,n}(s)$. 
Our criterion function will be
$$M_n(\t)= \|G_n(\t,\hat h_\t)\|$$
for some semi-norm $\|\cdot\|$ as described in the following assumption. 
\begin{enumerate}[label=(\textbf{N\arabic{*}})]
\item \label{N1} $\|\cdot\|$  is a semi-norm such that $\displaystyle\|\Gamma\|\leq c \sup_{y\in C\atop s\in[0,1]}|\Gamma(y,s)|$
for some constant $c>0$ and some compact set $C=[c_1,c_2]\subset\mathbb{R}$ with $c_1,c_2>0$ and $0\in C$, 
for all measurable functions $\Gamma:\mathbb{R}\times [0,1]\to\mathbb{R}$.
\end{enumerate}

For instance one can consider one of the following semi-norms, 
\begin{itemize}
\item[(i)] $\displaystyle\|\Gamma(y,s)\|=\sup_{s\in [0,1]\atop y\in C} |\Gamma(y,s)|$ 
\item[(ii)] $\displaystyle\|\Gamma(y,s)\|=\left(\int \Gamma(y,s)^2w(y,s)\,d(y,s)\right)^{1/2}$ for some integrable weight function $w:\mathbb{R}\times [0,1]\to\mathbb{R}_0^+$ with support included in $C\times [0,1]$
\item[(iii)] $\displaystyle\|\Gamma(y,s)\|=\sup_{s\in [0,1]}\left(\int \Gamma(y,s)^2w(y)\,dy\right)^{1/2}$ for some integrable weight function $w:\mathbb{R}\to\mathbb{R}_0^+$ with support included in $C$
\item[(iv)] $\displaystyle\|\Gamma(y,s)\|=\sup_{y\in C}\left(\int \Gamma(y,s)^2w(s)\,ds\right)^{1/2}$ for some integrable weight function $w:[0,1]\to\mathbb{R}_0^+$.
\end{itemize}
The first two semi-norms correspond to Kolmogorov-Smirnov and CramÚr-von Mises distances, respectively, while the last two are mixtures of both. 

\bigskip

Now we define the estimator $\hat \t$ of $\t_0$ as the minimizer of $M_n(\t)$ over $\Theta$, i.e.\
\begin{equation}\label{deftheta}
\hat\t=\arg\min_{\t\in\Theta}M_n(\t).
\end{equation}
For the following theory also the weaker condition $M_n(\hat\vartheta)\leq \inf_{\theta\in\Theta}M_n(\vartheta)+o_P(1)$
is sufficient. 
Note that 
\beq
\mathbb{P}\left(\Lt(Y_{i})-h(X_i)\leq y\mid X_i=x\right)&=&\mathbb{P}\left(\Lambda_0(Y_{i})\leq \Lambda_0(\Lt^{-1}(y+h(x)))\mid X_i=x\right)\\
&=&F_0\left(\Lambda_0(\Lt^{-1}(y+h(x)))-h_0(x)\right),
\eeq
which reduces to $F_0(y)$ for $\vartheta=\t_0$ and $h=h_0$. Now considering expectations we define 
\begin{eqnarray}\label{G-random}
G(\t,h)(y,s)&=&\int F_0\left(\Lambda_0(\Lt^{-1}(y+h(x)))-h_0(x)\right) I\{x\leq s\}f_X(x)\,dx\\
&&{}-\int F_0\left(\Lambda_0(\Lt^{-1}(y+h(x)))-h_0(x)\right)f_X(x)\, dx \, F_X(s)
\nonumber
\end{eqnarray}
 as deterministic counterpart of $G_n(\t,h)$. 
Further set $$M(\t)=\|G(\t, h_\t)\|$$ and note that $M(\t_0)=\|G(\t_0,h_{0})\|= 0$.
In section 4 we formulate assumptions under which $\hat\vartheta$ consistently estimates $\vartheta_0$. 

\subsection{The fixed design case}

In the fixed design model (\ref{model0}) we define the estimator for the boundary curve $h_0$ as
\[\tilde h_0(x)=\max\{\Lambda_0(Y_{i,n})|i=1,\ldots,n\text{ with }|\tin-x|\leq b_n\}\]
and obtain uniform consistency under the following modified assumptions. We set $x_{0,n}=0$ and $x_{n+1,n}=1$. 

\begin{enumerate}[label=(\textbf{A\arabic{*}'})]
\item \label{F1'} Model (\ref{model0}) holds with independent $\eps_{1,n},\dots,\eps_{n,n}$ with cdf $F_0$ ($\forall n$) such that $F_0(0)=1$, $F_0(-\Delta)<1$ for all $\Delta>0$.  
\item \label{X1'} The design points  $0<x_{1,n}<\dots <x_{n,n}<1$ are deterministic.
\setcounter{enumi}{3}
\item \label{B1'} Let $(b_n)_{n \geq 0}$ be a sequence of positive bandwidths that satisfies
$\lim_{n\to\infty}b_n=0$ and $\lim_{n \to \infty}\bar{\Delta}_n\log(n)/b_n=0$ for $\bar{\Delta}_n:=\max_{1 \leq i \leq n+1} \left(x_{i,n}-x_{i-1,n}\right)$. 
\end{enumerate}

\begin{lemma}\label{consistency} Under model  (\ref{model0}) with assumptions \ref{F1'}, \ref{X1'}, \ref{H1} and \ref{B1'} we have $\|\tilde h_0-h_0\|_{\infty}=o_P(1)$. 
\end{lemma}

The proof is given in section \ref{proof-consistency} of the online supplementary material. For general $\vartheta\in\Theta$ we define a consistent boundary curve  estimator  as
\[\tilde h_{\vartheta}(x)=\max\{\Lambda_{\vartheta}(Y_{i,n})|i=1,\ldots,n\text{ with }|\tin-x|\leq b_n\} = \Lt(\Lambda_0^{-1}(\tilde h_0(x))).\]
In analogy to (\ref{G_n-random}) we define, for any function $h:[0,1]\to\mathbb{R}$,
\begin{eqnarray}\label{G_n}
 G_n(\t, h)(y,s)&=& \frac 1n\sum_{i=1}^n I\{\Lt(Y_{i,n})-h(\tin)\leq y\} \big( I\{\tin\leq s\}-\hat F_{X,n}(s)\big),
\end{eqnarray}
where 
$$\hat F_{X,n}(s) = \frac 1n \sum_{i=1}^n I\{x_{i,n}\leq s\}.$$
The criterion function is again $M_n(\vartheta)=\|G_n(\vartheta,\hat h_\vartheta)\|$ where the smooth estimator $\hat h_\vartheta$ is defined accordingly as in \eqref{hatK} and with this the transformation parameter estimator is similar to (\ref{deftheta}). In order to consider the same deterministic $G$ as in (\ref{G-random}) an additional assumption is needed. 

\begin{enumerate}[label=(\textbf{A\arabic{*}''})]
\setcounter{enumi}{1}
\item \label{X1''} The design points  $0<x_{1,n}<\dots <x_{n,n}<1$ are deterministic. There exists a cdf $F_X$ with continuous density function $f_X:[0,1]\to \mathbb{R}$ which is bounded away from zero such that
$$\max_{i=1,\dots,n+1}\left|\int_{x_{i-1,n}}^{x_{i,n}}f_X(x)\,dx-{\frac 1n}\right|=o\left({\frac 1n}\right).$$
\end{enumerate}
Assumption \ref{X1''} is common in the literature on fixed design regression models. It allows the application of the mean value theorem for integrals to obtain, for some $\xi_{i,n}\in [x_{i-1,n},x_{i,n}]$, 
$$f_X(\xi_{i,n})(x_{i,n}-x_{i-1,n})=\int_{x_{i-1,n}}^{x_{i,n}}f_X(x)\,dx=\frac{1}{n}+o\left({\frac 1n}\right)$$
uniformly in $i=1,\dots,n$. Thus it follows from \ref{X1''} that $\bar{\Delta}_n$ in assumption \ref{B1'} has the exact rate $n^{-1}$ and therefore assumption \ref{B1'} reduces to \ref{B1}.  Further the following Riemann sum approximations for bounded integrable functions $\varphi$ can be applied to get
\begin{eqnarray}\nonumber
\frac 1n \sum_{i=1}^n \varphi(x_{i,n}) &=& \sum_{i=1}^n \varphi(x_{i,n}) (f_X(\xi_{i,n})(x_{i,n}-x_{i-1,n})+o({\textstyle{\frac 1n}})) \\
& =&\int \varphi(x)f_X(x)\,dx+o(1).  \label{riemann}
\end{eqnarray}
In the next section we state conditions under which $\hat\vartheta=\arg\min_{\vartheta\in\Theta}M_n(\vartheta)$ consistently estimates $\vartheta_0$.

\section{Main result}

To prove consistency of the estimator for the transformation parameter we need the following additional assumptions. Please note that assumption \ref{M1} implies identifiability of the transformation $\Lambda_0$ in the class $\mathcal{L}$.

\begin{enumerate}[label=(\textbf{B\arabic{*}})]
\item \label{M1} For every $\delta>0$ there exists some $\epsilon>0$ such that $\inf_{\|\t-\t_0\|>\delta} M(\t)\geq\epsilon$.
\item \label{L0} $\mathcal{L}=\{\Lambda_\t\mid \t\in\Theta\}$ is a class of strictly increasing continuous functions $\mathbb{R}\to\mathbb{R}$. 
\item \label{L1} Let $S=\{\Lambda_0^{-1}(h_0(x))\mid x\in[0,1]\}$. Then the class $\mathcal{L}_S=\{\Lambda_\t|_S\mid \t\in\Theta\}$ is pointwise bounded and uniformly equicontinuous, i.e.\ $\sup_{\t\in\Theta}|\Lambda_\t(y)|<\infty$ for all $y\in S$, and  for  every $\epsilon>0$ there exists some $\delta>0$ such that $\sup_{\t\in\Theta}|\Lambda_\t(y)-\Lambda_\t (z)|<\epsilon$ for all $y,z\in S$ with $|y-z|\leq\delta$. 
\item \label{L2} The class $\mathcal{L}_{\tilde S}^1=\{\Lambda_0\circ\Lambda_\t^{-1}|_{\tilde S}\mid \t\in\Theta\}$
is pointwise bounded and uniformly equicontinuous for $\tilde S=\{z+h_\t(x)\mid z\in C_\tau, \t\in\Theta, x\in[0,1]\}$ with $C_\tau=[c_1-\tau,c_2+\tau]$ (for $C=[c_1,c_2]$ from \ref{N1}) for some $\tau>0$, i.e.\  
 $\sup_{\t\in\Theta}|\Lambda_0(\Lambda_\t^{-1}(z))|<\infty$ for all $z\in \tilde S$, and 
for every $\delta>0$ there exists some $\gamma>0$ such that 
$\sup_{\t\in\Theta}|\Lambda_0(\Lambda_\t^{-1}(x))-\Lambda_0(\Lambda_\t^{-1}(z))|\leq\delta$ for all $x,z\in \tilde S$ with $|x-z|\leq\gamma$. 
\item \label{F2} For some $\tau>0$, $F_0$ is uniformly continuous on the set $\tilde C=\{\Lambda_0(\Lambda_\t^{-1}(y+a+h_\t(x)))-h_0(x)\mid y\in C, \t\in\Theta,x\in[0,1],|a|\leq \tau\}$ (with $C$ from \ref{N1}), i.e.\ for every $\epsilon>0$ there is some $\delta>0$ such that $|F_0(y)-F_0(z)|<\epsilon$ if $|y-z|\leq\delta$, $y,z\in \tilde C$. 
\item\label{smoothing} $K$ is a density with support $[-1,1]$ and $b_n\searrow 0$, $nb_n\to\infty$. 
 \end{enumerate}
 Let us now make few comments regarding these assumptions. 
 
 \begin{itemize}
 	\item  \ref{M1} is a common assumption in M-estimation and needed for uniqueness of the true parameter. 
\item \ref{L0} implies the existence of continuous inverse functions $\Lambda_\t^{-1}$. 
Further note that uniform equicontinuity and pointwise boundedness imply totally boundedness by the ArzelÓ-Ascoli theorem. Thus for each $\epsilon$ there is a finite  covering of the classes $\mathcal{L}_S$ from \ref{L1} and $\mathcal{L}_{\tilde S}^1$  from \ref{L2} with balls of radius $\epsilon$ with respect to the sup norm. Thus also   the sup norm bracketing numbers of those classes are finite, i.e.
\begin{eqnarray}\label{bracketing}
N_{[\,]}(\epsilon,\mathcal{L}_S,\|\cdot\|_\infty)<\infty, \quad N_{[\,]}(\epsilon,\mathcal{L}^1_{\tilde S},\|\cdot\|_\infty)<\infty  \mbox{ for all } \epsilon>0
\end{eqnarray}
(see, e.g., Lemma 9.21 in Kosorok, 2008).

\item  \ref{L1}--\ref{F2} can be seen as minimal assumptions on the class $\mathcal{L}=\{\Lambda_\vartheta\mid\vartheta\in\Theta\}$ and $F_0$. As typically the sets $S$, $\tilde S$ and $\tilde C$ are unknown, the assumptions can be replaced by stronger assumptions that hold on all compact sets.
Besides, working on compact set transformation parameter, assumptions \ref{L1}--\ref{L2} hold for most of transformations used in practice such as the Box and Cox transformations (see Box and Cox, 1964) (suitably modified taking into account the data range), the exponential transformations (see Manly, 1976), the sinh-arcsinh transformations  (see Jones and Pewsey, 2009). For instance, with regard to Yeo-Johnson transformations, when $\vartheta \in \Theta=[0,2]$, $\Lambda_{\vartheta}: \mathbb{R}\to \mathbb{R}$ defines a bijective map (see Remark \ref{yeo-johnson}) and both $\Lambda_\t$ and $\Lambda_\t^{-1}$ have uniform bounded derivatives on compact sets so that one may show that they fulfill assumptions \ref{L1}--\ref{L2} using the mean value theorem to $\Lambda_\t$ and $\Lambda_\t^{-1}$. 
Further under stronger assumptions on the smoothness of $F_0$, $\Lambda_\t$ and $\Lambda_0\circ\Lambda_\t^{-1}$, the theoretical results can be generalized to semi-norms that are not restricted to a compact $C\times[0,1]$ as in assumption \ref{N1}. 
\item \ref{smoothing} is standard in kernel smoothing and is needed for the smoothed estimator $\hat h_\vartheta$ to be consistent. While we noticed in the simulations that slight smoothing improves the procedure, the following theorem still holds when $\hat h_\vartheta$ is replaced by the non-smooth estimator $\tilde h_\vartheta$. Assumption \ref{F2} holds, e.g.\ for H÷lder-continuous distribution functions $F_0$. 

\end{itemize}

The following theorem states consistency of the transformation parameter estimator. 

\begin{theo}\label{theo}
{\bf (i).}
(The random design case.) 
Assume  model (\ref{model0-random}) under assumptions \ref{F1}--\ref{B1}, \ref{N1}, \ref{M1}--\ref{smoothing}. Then $\hat\t$ is a consistent estimator, i.e.\ $\hat\t-\t_0=o_P(1)$.

{\bf (ii).}
(The fixed design case.) 
Assume  model (\ref{model0})  under assumptions \ref{F1'}, \ref{X1''}, \ref{H1}, \ref{B1}, \ref{N1}, \ref{M1}--\ref{smoothing}. Then $\hat\t$ is a consistent estimator, i.e.\ $\hat\t-\t_0=o_P(1)$.
\end{theo}

The proof for the random design case is given in section \ref{proof-theo-random} of the appendix and the proof for the fixed design case in section \ref{proof-theo-fixed} of the supplement. 
One basic ingredient is the following result, which is proven in section \ref{proof-lemma-smoothing} of the appendix for the random design case. The proof for the fixed design case is analogous. 

\begin{lemma}\label{lemma-smoothing}
{\bf (i).}
(The random design case.) 
Under model  (\ref{model0-random}) with assumptions \ref{F1}--\ref{B1}, \ref{L0}, \ref{L1}, \ref{smoothing}, we have $\sup_{\vartheta\in\Theta}\|\hat h_\vartheta-h_\vartheta\|_\infty=o_P(1)$.

{\bf (ii).}
(The fixed design case.) 
Under model  (\ref{model0}) with assumptions \ref{F1'}, \ref{X1'}, \ref{H1}, \ref{B1'}, \ref{L0}, \ref{L1}, \ref{smoothing}, we have  $\sup_{\vartheta\in\Theta}\|\hat h_\vartheta-h_\vartheta\|_\infty=o_P(1)$.
\end{lemma}
\bigskip

The consistency result in Theorem \ref{theo} should be seen as a first step in the analysis of transformation boundary regression models. An interesting and challenging topic for future research is to derive an asymptotic distribution of $\hat\vartheta-\vartheta_0$ (properly scaled) and to investigate the asymptotic influence of the estimation on subsequent procedures based on the transformed data. This is beyond the scope of the paper as yet there are no results on the uniform asymptotic distribution of $\tilde h_0-h_0$ in the literature.\newline 

We finally highlight that under the further condition  \ref{X3bis} defined below regarding the regularity of the boundary curve, we obtain as a corollary of Theorem \ref{theo} the consistency of the estimator $\hat h_{\hat \vartheta}$ of the boundary curve. 

\begin{enumerate}[label=(\textbf{A\arabic{*}'})]
	\setcounter{enumi}{2}
	\item \label{X3bis}  $\vartheta_0$ is an inner point of a convex parameter space $\Theta$ and $h_\vartheta$ is continuously differentiable with respect to $\vartheta$. Besides, we assume that there exists some $\delta>0$ such that $$\sup_{x \in [0,1]} \sup_{\Vert \vartheta - \vartheta_0 \Vert<\delta}  \left\| \frac{\partial h_\vartheta (x)}{\partial \vartheta} \right\|< \infty. $$
\end{enumerate}
\begin{cor}
	{\bf (i).}
	(The random design case.) 
	Assume  model (\ref{model0-random})  holds under assumptions \ref{F1}, \ref{X1}, \ref{X3bis}, \ref{B1}, \ref{N1} and \ref{M1}--\ref{smoothing}. Then $\hat h_{\hat \vartheta}$ is a consistent estimator of $ h_{ \vartheta_0}$, i.e.\ $\Vert \hat h_{\hat \vartheta} - h_{ \vartheta_0} \Vert_{\infty} =o_P(1)$.
	
	{\bf (ii).}
	(The fixed design case.) 
	Assume  model (\ref{model0}) holds under assumptions \ref{F1'}, \ref{X1''}, \ref{X3bis}, \ref{B1}, \ref{N1} and \ref{M1}--\ref{smoothing}. Then $\hat h_{\hat \vartheta}$ is a consistent estimator of $\ h_{ \vartheta_0}$, i.e.\ $\Vert \hat h_{\hat \vartheta} - h_{ \vartheta_0} \Vert_{\infty} =o_P(1)$.
\end{cor}

{\bf Proof.} We only prove {\bf (i)} since the proof of {\bf (ii)} is identical. 
	Observe first that 
	\begin{align*}
\Vert \hat h_{\hat \vartheta} - h_{ \vartheta_0} \Vert_{\infty} \leq \Vert \hat h_{\hat \vartheta} - h_{ \hat  \vartheta} \Vert_{\infty} + \Vert h_{ \hat \vartheta} - h_{ \vartheta_0} \Vert_{\infty}.
\end{align*}
The first term in the right hand side of the above inequality goes to $0$ in probability from Lemma \ref{lemma-smoothing} since the consistency holds uniformly over $\vartheta \in \Theta$. Regarding the second term, applying the mean value theorem, there exists some $\vartheta^*(x)$ on the line between $\hat \vartheta$ and $\vartheta_0$ such that

$$ \Vert  h_{\hat \vartheta} - h_{  \vartheta_0} \Vert_{\infty} = \sup_{x \in [0,1]} \left| \frac{\partial h_\vartheta (x)^T}{ \partial \vartheta} |_{\vartheta=\vartheta^*(x)} (\hat \vartheta-\vartheta_0) \right|.$$
From Theorem \ref{theo},  $\hat \vartheta-\vartheta_0=o_P(1)$ which concludes the proof under the assumption \ref{X3bis}.
\hfill $\Box$


\section{Simulations}

To study the small sample behavior, we generate data as $Y=\Lambda_{\vartheta_0}^{-1}(h_0(x)+\varepsilon)$ using the Yeo-Johnson transformation for different values of $\vartheta_0$. We focus on the equidistant design framework and examine the two regression functions $h_0(x)= 10(x-\textstyle{\frac 12})^2$ and $h_0(x)= \textstyle{\frac 12}\sin(2\pi x)+4x$ for two different error distributions, namely the Weibull distribution with scale parameter $1$ and shape parameter $3$ and the exponential distribution with mean $1/3$.  We consider samples of size $n=50$ and $n=100$. It means that we investigate the following four models  
\begin{eqnarray}
h_0(x)= 10(x-\textstyle{\frac 12})^2   && \text{with} \ \ \ \varepsilon\sim \mbox{Weibull}(1,3)  \label{simu1}\\
h_0(x)= 10(x-\textstyle{\frac 12})^2  && \text{with} \ \ \ \varepsilon\sim \mbox{Exp}(3)\label{simu2}\\
h_0(x)= \textstyle{\frac 12}\sin(2\pi x)+4x && \text{with} \ \ \ \varepsilon\sim \mbox{Weibull}(1,3)  \label{simu3}\\
h_0(x)=  \textstyle{\frac 12}\sin(2\pi x)+4x && \text{with} \ \ \ \varepsilon\sim \mbox{Exp}(3). \label{simu4}
\end{eqnarray}
Figures \ref{graphic1} and \ref{graphic2} show realizations of models (\ref{simu1}) and (\ref{simu3}).
The bandwidth $b_n=n^{-1/3}$ is chosen accordingly to Drees et al.\ (2018) and simulations are based on $1000$ iterations. We use the Epanechnikov kernel to smooth the boundary curve estimator and compare the results for two smoothing parameters $a_n=b_n/2$ and $a_n=b_n/20$.  
The transformation parameter estimator is as in (\ref{deftheta}) on the interval $[-0.5,2.5]$, where the semi-norm in the criterion function $M_n(\vartheta)$ is chosen as in $(i)$, $(ii)$, $(iii),$ and $(iv)$ in the examples of Condition \ref{N1}. In the following we denote the according estimators as  TKS, TCM, TKSCM and TCMKS. Here, TKS and TCM refer to Kolmogorov-Smirnov and CramÚr-von Mises distances respectively, while TKSCM and TCMKS are mixtures of both.
For simplicity, the weight functions are chosen identically equal to $1$ in all the settings, i.e., $w(y,s) =1$ for all $ (y,s) \in \mathbb{R} \times [0,1] $, $w(y)=1$ for all $y \in \mathbb{R}$ and $w(s)=1$ for all $s \in [0,1]$ in $(ii)$, $(iii)$ and $(iv)$, respectively (although for the theory we assumed a compact support). \newline  
\noindent We sum up the simulation results in the following 8 tables. Tables \ref{Table1}, \ref{Table7} and Tables \ref{Table2}, \ref{Table8} deal with Models \eqref{simu1} and \eqref{simu4} for $a_n=b_n/2$ and $a_n=b_n/20$, respectively, whereas Tables  \ref{Table3}, \ref{Table5} and Tables \ref{Table4}, \ref{Table6} in the supplement show the results for Models \eqref{simu2} and \eqref{simu3}.  In Figure \ref{est}, we have represented the density function of each estimator for the Model \eqref{simu1} when $\t_0=0.5$ with $n=100$ and $a_n=b_n/20$, which corresponds to the settings of Table \ref{Table2}. To assess the performance of our estimates, we provide for each estimator the mean, the median and the Mean Integrated Squared Error (MISE) in brackets for five values of the true parameter $\t_0=0,0.5,1,1.5,2$. The best-performing one regarding the mean (respectively
the MISE) is highlighted in bold (respectively is underlined).  

Looking at the MISE, it turns out that the estimator using the Cram\'{e}r-von Mises distance (TCM) out-performs in many cases even when it does not out-perform the mean; see Tables \ref{Table1} and \ref{Table2} when $n=100$ for instance. Besides, as it is intented, results are better in most of the cases when the sample size $n$ increases. However, this does not hold for every case. For instance, one may see in Table \ref{Table2} that for the second estimator TCM, most of the results are better for $n=50$ than for $n=100$. 
This might relate to a sensitivity with respect to the choices of bandwidth and smoothing parameter. 
 A lot of criteria may be used to judge the performance of the estimators. We deal here with the mean, the median and the MISE but we emphasize that using different criteria (e.g.\ median absolute deviation, mode or even graphical analysis)  could give different results concerning the comparison of the methods. For instance, results in Table \ref{Table7} for $n=100$, $a_n=b_n/2$ and $\t_0=1$ are quite not accurated regarding the mean (e.g. $0.845$ for the TCM). Nevertheless, looking at Figure \ref{est2}, it appears that the plots of the densities look satisfactory.


\noindent It is clear that the TCM and the TKSCM out-perform in the Model \eqref{simu3} and in the Model \eqref{simu4}, respectively. Nonetheless, in a general setting, we are not able to state which estimator performs better since it depends first on the criterium selected to judge the performance but more importantly on the choice of the bandwidths and the smoothing parameter. 

Finally, we recall that the aim of this work is to reduce the dependence between the covariates and the errors. As one can see in Table \ref{TableCor2paper}, although the estimation of $\vartheta_0$ is less good than expected in the model \eqref{simu4}, the correlations between the covariates and the errors (after transformation) are very small; see also Table \ref{TableCor1paper} in the supplement for the correlations in the model \eqref{simu3}. We obtain similar results for the random design case.

\begin{appendix}

\section{Proofs of asymptotic results in the random covariate case}

For the proofs of the asymptotic results let us fix some notation: $\lfloor \cdot \rfloor$ and $\lceil \cdot \rceil$ are the floor and ceiling functions respectively; $\bar{F}=1-F$ denotes the survival function associated to a cdf $F$; $X_1 \overset{d}{=} X_2$ means that two random variables $X_1,X_2$ share the same distribution; $a_n \underset{n \to \infty}{\sim} b_n$ holds if there exists a constant $ c>0$ such that $\lim_{n \to \infty}a_n/b_n=c$ for two sequences $(a_n)_{n \geq 1}$ and $(b_n)_{n \geq 1}$ of nonnegative numbers; $A^c$ is the complement of a set $A$.\newline 

\noindent In the following we give the proofs of our results in the random design case whereas the proofs for the fixed design case can be found in the online supplementary material.


\subsection{Proof of Lemma \ref{consistency-random}}\label{proof-consistency-random}
At first we need the following intermediary lemma. 

\begin{lemma}\label{lem:model0-random}
	Assume model (\ref{model0-random}) holds with assumptions \ref{F1}, \ref{X1} and \ref{B1}. Then, we have
	\begin{eqnarray*}\label{Lemma2dot2}
		\sup_{x \in [0,1]}  \min_{ \substack{ i \in \{1,\ldots,n\} \\ |X_{i}-x|\leq b_n }} |\varepsilon_{i}|=o_P(1).
	\end{eqnarray*}
\end{lemma}

\noindent {\bf Proof.}	
	For $n \geq 1$ denote $X_{(1)}<X_{(2)} < \cdots < X_{(n)} $ the order statistics of the random design sample $X_1, X_2,\ldots, X_n$. 
Let $\pi$ be the random permutation of $\{ 1, \ldots, n\}$ such that $X_{(i)}=X_{\pi(i)}, i=1, \ldots,n$. Due to the independence between the errors and the covariates under 
\ref{F1}, $\eps_{\pi(1)},\dots,\eps_{\pi(n)}$ are iid with cdf $F_0$. Let $Z_i=-\eps_{\pi(i)}$, $i=1,\dots,n$, then $Z_1,\dots,Z_n$ are iid with cdf $U$ with $U(x)=1-F_0(-x)$ and we need to show that
	\begin{align}\label{eq1_lem_RD}
	\lim_{n \to \infty} \mathbb{P}\left( \sup_{x \in [0,1]}  \min_{ \substack{ i \in \{1,\ldots,n\} \\ |X_{(i)}-x|\leq b_n }} Z_{i}> \epsilon \right)=0, \ \ \ \epsilon>0.
	\end{align}	
Define for $n \geq 1$ the event 
$$\Omega_n=\left\{ \inf_{x\in [0,1]}\sum_{i=1}^n I\{|X_i-x|\leq b_n\}\geq C nb_n  \right\}$$
for a suitable constant $C>0$ specified later. Note that on $\Omega_n$, there are at least $Cnb_n$ covariates in each of the intervals $[x-b_n,x+b_n]$. We will first show that $\lim_{n\to\infty}\mathbb{P}(\Omega_n)=1$. To this end, for $n \geq 1$ let $f_{n,x}(z)=I\{|x-z|\leq b_n\}$ and note that
\begin{eqnarray}
\nonumber
	\inf_{x\in [0,1]}\frac{1}{n}\sum_{i=1}^n I\{|X_i-x|\leq b_n\}&\geq& \inf_{x\in [0,1]} \mathbb{P}(|X_1-x|\leq b_n)\\
	&&{}-\sup_{x\in[0,1]}\Big|\frac 1n\sum_{i=1}^n (f_{n,x}(X_i)-\mathbb{E}[f_{n,x}(X_i)])\Big|.
	\label{sternchen1}
\end{eqnarray}
Applying the mean value theorem of integration, it follows that
\begin{equation}\label{milton}
2b_n\sup_{x\in[0,1]}f_X(x)\geq \mathbb{P}(|X_1-x|\leq b_n) = \int_{\max(0,x-b_n)}^{\min(1,x+b_n)}f_X(x)\,dx \geq b_n \inf_{x\in[0,1]} f_X(x) .
\end{equation}
Then, there exists a constant $C_1>0$, which actually corresponds to the lower bound of the density function $f_X$ involved in Assumption \ref{X1} such that 
\begin{equation}\label{sternchen2}
\mathbb{P}(|X_1-x|\leq b_n) \geq C_1b_n,
\end{equation}
uniformly over $x \in [0,1].$

\noindent Fix $n \geq 1$ and denote $P_nf_{n,x}:=\frac{1}{n}\sum_{i=1}^n f_{n,x}(X_i)$ and $Pf_{n,x}:=\mathbb{E}[f_{n,x}(X_1) ]$ so that $P_n$ and $P$ refer to the empirical measure and the distribution of the random design sample $X_1, \ldots, X_n$, respectively.
By (\ref{milton})
$Pf^2_{n,x} = \mathbb{E}[I\{|X-x|\leq b_n \} ]\leq 2 C_2 b_n$,
where $C_2:=\sup_{x \in [0,1]} f_X(x)$, which is finite under  \ref{X1}. Moreover, since $|f_{n,x}(X)|\leq 1$ and the assumption on the covering number is fulfilled (see Example 38 and Problem 28 to be convinced in Pollard (1984)), Theorem 37 in Pollard (1984, p. 34) holds and we have
\begin{eqnarray*}
	\sup_{x\in[0,1]}\Big|\frac 1n\sum_{i=1}^n (f_{n,x}(X_i)-\mathbb{E}[f_{n,x}(X_i)])\Big|=o(b_n).
\end{eqnarray*}
From this together with (\ref{sternchen1}) and (\ref{sternchen2}) it follows that $\lim_{n\to\infty}\mathbb{P}(\Omega_n)=1$. It means that for any sub-interval $I_n:=[x-b_n,x+b_n]$, there are at least $Cnb_n$ random design points with probability converging to $1$. 

Then, for all $y>0$, we have with $d_n:=\lceil Cnb_n \rceil $
\begin{eqnarray*}
	\mathbb{P}\left( \sup_{x \in [0,1]}  \min_{ \substack{ i \in \{1,\ldots,n\} \\ |X_{(i)}-x|\leq b_n }} Z_{i}> y \right  )  &\leq&   \mathbb{P}\left( \left\lbrace  \max_{j \in \{1,\ldots,n-d_n\}} \min_{i \in \{j,\ldots,j+d_n\}} Z_{i}> y \right\rbrace
 \cap \Omega_n	\right) + \mathbb{P}\left( \Omega_n^c \right)\\
	&\leq& \sum_{j=1}^{n-d_n}\mathbb{P}\left(\min_{i\in\{j,\ldots,j+d_n\}}Z_i>y\right)+ \mathbb{P}\left( \Omega_n^c \right)\\
	&=& (n-d_n)\mathbb{P}\left(\min_{i\in\{1,\ldots,d_n+1\}}Z_i>y\right)+ \mathbb{P}\left( \Omega_n^c \right)\\
	&=&(n-d_n)\overline{U}(y)^{d_n+1}+ \mathbb{P}\left( \Omega_n^c \right).
	\end{eqnarray*}
Thus it remains to show that for all $\epsilon>0$
\[(n-d_n)\overline{U}(\epsilon)^{d_n+1}\nto 0\]
which is true since  $d_n \underset{n \to \infty}{\sim} nb_n$ and
\beq
nb_n\log(\overline{U}(\epsilon))+\log(n-d_n)&\leq&nb_n\log(\overline{U}(\epsilon))+\log(n)\\
&=&\log(n)\left(\frac{nb_n}{\log(n)}\log(\overline{U}(\epsilon))+1\right)\\
&\nto& -\infty
\eeq
since $\overline{U}(\epsilon)<1$ under \ref{F1} and $\frac{nb_n}{\log(n)} \nto \infty$ under \ref{B1}. This concludes the proof.
\hfill $\Box$

\bigskip

We are now ready to prove Lemma \ref{consistency-random}.

\noindent
{\bf Proof of Lemma \ref{consistency-random}.}
	On the one hand, we have
	\begin{eqnarray}\label{eq1_th1}
	\sup_{x \in [0,1]} \left(\tilde{h}_0(x)-h_0(x) \right) \nonumber & = & \sup_{x \in [0,1]} \left( \max_{ \substack{ i \in \{1,\ldots,n\} \\ |X_{i}-x|\leq b_n }} \left\{ h_0(X_{i})+\varepsilon_{i} - h_0(x)\right\}\right)\\ 
	& \leq & \sup_{x \in [0,1]} \left( \max_{ \substack{ i \in \{1,\ldots,n\} \\ |X_{i}-x|\leq b_n }} \left\{ h_0(X_{i}) - h_0(x)\right\}\right)\\ \nonumber
	& \leq & \sup_{|t-x| \leq b_n} |h_0(t)-h_0(x)|\\ \nonumber
	& = &o(1),
	\end{eqnarray}
	since the errors $(\varepsilon_{i})_{1 \leq i \leq n}$ are nonpositive and $h_0$ is continuous on the compact set $[0,1]$ and thereby uniformly continuous under \ref{H1}. On the other hand, 
	\begin{eqnarray}\label{eq2_th1}
	\sup_{x \in [0,1]} \left(h_0(x) -\tilde{h}_0(x)\right) \nonumber & =& \sup_{x \in [0,1]} \left( h_0(x) -\max_{ \substack{ i \in \{1,\ldots,n\} \\ |X_{i}-x|\leq b_n }}\left\{ h_0(X_{i})+\varepsilon_{i,n}\right\}\right)\\ 
	& = &\sup_{x \in [0,1]} \left( \min_{ \substack{ i \in \{1,\ldots,n\} \\ |X_{i}-x|\leq b_n }} \left\{ h_0(x) - h_0(X_{i})-\varepsilon_{i,n}\right\}\right)\\ \nonumber
	& \leq & \sup_{x \in [0,1]} \left( \min_{ \substack{ i \in \{1,\ldots,n\} \\ |X_{i}-x|\leq b_n }} \left\{ -\varepsilon_{i} \right\}\right) + \sup_{|t-x| \leq b_n} |h_0(t)-h_0(x)|\\ \nonumber
	& = & o_P(1),
	\end{eqnarray}
	with Lemma \ref{lem:model0-random}.
		Finally, combining equations \eqref{eq1_th1} and \eqref{eq2_th1}, it follows that
	\begin{eqnarray*}
	\Vert h_0-\tilde{h}_0\Vert_{\infty} &=& \sup_{x \in [0,1]} \left|\tilde{h}_0(x)-h_0(x) \right|\\
	& =& \sup_{x \in [0,1]} \left( \max \left\{ \tilde{h}_0(x)-h_0(x), h_0(x)-\tilde{h}_0(x) \right\}  \right)\\
	& \leq & \max \left\{  \sup_{x \in [0,1]} (\tilde{h}_0(x)-h_0(x)), \sup_{x \in [0,1]} (h_0(x)-\tilde{h}_0(x)) \right\}\;=\;o_P(1),
	\end{eqnarray*}
	which is the desired result.

 \hfill $\Box$

\subsection{Proof of Lemma \ref{lemma-smoothing}}\label{proof-lemma-smoothing}

Let $\epsilon>0$. 
Note that $\Lambda_0^{-1}\circ h_0$ is uniformly continuous due to assumptions \ref{H1} and \ref{L0}. Thus with $h_\vartheta=\Lambda_\vartheta\circ\Lambda_0^{-1}\circ h_0$ and assumption \ref{L1} it follows that there exists some $\delta>0$ such that $\sup_{\vartheta\in\Theta}|h_\vartheta(x)-h_\vartheta(y)|\leq\epsilon$ if $|x-y|\leq\delta$. Now let $n$ be large enough such that $b_n\leq\delta$. Then due to the definition of $\hat h_\vartheta$ and $\mbox{supp}(K)=[-1,1]$ one obtains
$$\|\hat h_\vartheta-h_\vartheta\|_\infty\leq \|\tilde h_\vartheta-h_\vartheta\|_\infty+\sup_{x\in[0,1]}\Big|\frac{\sum_{i=1}^n ( h_{\vartheta}(X_i)- h_{\vartheta}(x))K(\frac{x-X_i}{b_n})}{\sum_{i=1}^n K(\frac{x-X_i}{b_n})}\Big|\leq \|\tilde h_\vartheta-h_\vartheta\|_\infty+\epsilon.$$  
From  Lemma \ref{consistency-random} we have  $\|\tilde h_0-h_0\|_{\infty}=o_P(1)$ and thus with assumption \ref{L1} it follows that
\begin{eqnarray*}
\sup_{\t\in\Theta}\| \tilde h_\t- h_\t\|_\infty&=&\sup_{\t\in\Theta}\sup_{x\in[0,1]}|\Lt( \Lambda_0^{-1}( \tilde h_0(x)))-\Lt( \Lambda_0^{-1}( h_0(x)))|
\;=\; o_P(1)
\end{eqnarray*}
and therefore the assertion of the lemma. \hfill $\Box$

\subsection{Proof of Theorem \ref{theo}}\label{proof-theo-random}

By the argmax theorem applied to the criterion function $M_n(\vartheta)$ multiplied by $(-1)$ and using assumption \ref{M1} it suffices to show that
$$\sup_{\vartheta\in\Theta}|M_n(\vartheta)-M(\vartheta)|=o_P(1)$$
(see\ Kosorok, 2008, Theorem 2.12(i)). 
To obtain this, note that
\begin{eqnarray*}
\sup_{\t\in\Theta}|M_n( \t)-M(\t)|&\leq &  \sup_{\t\in\Theta}\|G_n(\t,\hat h_{ \t})-\bar G_n(\t,\hat h_{\t})\| + \sup_{\t\in\Theta}\|\bar G_n(\t,\hat h_{\t})-G(\t,\hat h_{\t})\|\\
&&{}+
\sup_{\t\in\Theta} \|G(\t,\hat h_{\t})-G(\t, h_{\t})\|, 
\end{eqnarray*}
where
\begin{eqnarray}\label{G_n-bar}
\bar G_n(\t,h)(y,s)
&=& \frac 1n\sum_{i=1}^n I\{\Lt(Y_{i})-h(X_i)\leq y\} I\{X_i\leq s\}\\
&&{} -F_X(s)\frac 1n\sum_{i=1}^n I\{\Lt(Y_{i})-h(X_i)\leq y\}.\nonumber
\end{eqnarray}
Note that for any deterministic function $h$ we have $\mathbb{E}[\bar G_n(\t,h)]=G(\t,h)$.
The assertion of the theorem follows from 
$$\sup_{\t\in\Theta}\|G_n(\t,\hat h_{ \t})-\bar G_n(\t,\hat h_{\t})\|\leq \sup_{s\in[0,1]}|\hat F_{X,n}(s)-F_X(s)|=o_P(1)$$
and Lemmas \ref{lem2-random} and \ref{lem1-random}. \hfill $\Box$

\begin{lemma}\label{lem2-random} Under the assumptions of Theorem \ref{theo} (i),
$$\sup_{\t\in\Theta}\|\bar G_n(\t,\hat h_{ \t})- G(\t,\hat h_{ \t})\| =o_P(1).$$
\end{lemma}

\noindent
{\bf Proof.}
From Lemma \ref{lemma-smoothing} follows the existence of some deterministic sequence $a_n\searrow 0$ such that the probability of the event
\begin{equation}\label{E_n}
\sup_{\t\in\Theta}\| h_\t-\hat h_\t\|_\infty\leq a_n
\end{equation}
converges to one. Thus we assume in what follows that (\ref{E_n}) holds. 

We only consider the difference between the first sum in the definition of $\bar G_n(\t, h)$ (see (\ref{G_n-bar})) and the first integral in the definition of  $G(\t, h)$  (see (\ref{G-random})). The difference between the second sum and the second integral can be treated similarly. 
Applying (\ref{E_n}) the first sum in $\bar G_n(\t, \hat h_\t)(y,s)$ can be nested as 
\begin{eqnarray*}
&&\frac 1n\sum_{i=1}^n  I\{\Lt(Y_{i})-h_\t(X_i)\leq y-a_n\}I\{X_i\leq s\}\\
&\leq& \frac 1n\sum_{i=1}^n  I\{\Lt(Y_{i})-\hat h_\t(X_i)\leq y\}I\{X_i\leq s\}\\
&\leq& \frac 1n\sum_{i=1}^n  I\{\Lt(Y_{i})-h_\t(X_i)\leq y+a_n\}I\{X_i\leq s\}
\end{eqnarray*}
while the first integral in $  G(\t, \hat h_\t)(y,s)$ can be nested as 
\begin{eqnarray*}
&& \int F_0\left(\Lambda_0(\Lt^{-1}(y-a_n+h_\t(x)))-h_0(x)\right)I\{x\leq s\} f_X(x)\,dx\\
&\leq& \int F_0\left(\Lambda_0(\Lt^{-1}(y+\hat h_\t(x)))-h_0(x)\right)I\{x\leq s\} f_X(x)\,dx\\
&\leq& \int F_0\left(\Lambda_0(\Lt^{-1}(y+a_n+h_\t(x)))-h_0(x)\right)I\{x\leq s\} f_X(x)\,dx.
\end{eqnarray*}
Thus we have to consider
\begin{eqnarray*}
&&H_{n,\t}^{(1)}(y,s) \\
&=&\frac 1n\sum_{i=1}^n  \Big(I\{\Lt(Y_{i})-h_\t({X_i})\leq y+a_n\}I\{X_i\leq s\} \\
&&{}\qquad-  \int F_0\left(\Lambda_0(\Lt^{-1}(y+a_n+h_\t(x)))-h_0(x)\right)I\{x\leq s\} f_X(x)\,dx\Big)\\
\\
&&H_{n,\t}^{(2)}(y,s) \\
&=&\int  \Big(F_0\left(\Lambda_0(\Lt^{-1}(y+a_n+h_\t(x)))-h_0(x)\right)\\
&&{}\qquad-
F_0\left(\Lambda_0(\Lt^{-1}(y+h_\t(x)))-h_0(x)\right)\Big)I\{x\leq s\} f_X(x)\,dx\Big)
\end{eqnarray*}
and the same terms with $y+a_n$ replaced by $y-a_n$, which can be treated completely analogously. We have to show that
$\sup_{\t\in\Theta}\|H_{n,\t}^{(1)}\|=o_P(1)$ and $\sup_{\t\in\Theta}\|H_{n,\t}^{(2)}\|=o(1)$. 




Recall condition \ref{N1} and note that $\sup_{\t\in\Theta}\sup_{s\in [0,1]\atop y\in C}|H_{n,\t}^{(2)}(y,s)|=o(1)$ follows from uniform continuity of $F_0$ and of $\Lambda_0\circ \Lt^{-1}$ uniformly in $\t$  (see \ref{F2} and \ref{L2}), from the representation $h_\t=\Lambda_\t\circ\Lambda_0^{-1}\circ h_0$ and uniform continuity of $\Lambda_\t$ uniformly in $\t$ (see \ref{L1}), and $a_n\to 0$. 

Let $n$ be large enough such that $|a_n|\leq \tau$ for $\tau$ both from \ref{F2} and \ref{L2}. Now to prove $\sup_{\t\in\Theta}\|H_{n,\t}^{(1)}\|=o_P(1)$ note that
$$\sup_{\t\in\Theta}\|H_{n,\t}^{(1)}\|\leq \sup_{f\in\mathcal F}|P_nf-Pf|,$$
where $P_n$ denotes the empirical measure of $(X_1,Y_1),\dots,(X_n,Y_n)$, and $P$ the measure of $(X_1,Y_1)$, and 
$$\mathcal F=\{ (x,y)\mapsto I\{\Lt(y)-h_\t(x)\leq z\}I\{x\leq s\} \mid \t\in\Theta, s\in[0,1], z\in C_\tau\}$$
with $C_\tau$ as in assumption \ref{L2}. 
The assertion follows from the Glivenko-Cantelli theorem as stated in Theorem 2.4.1 in van der Vaart and Wellner (1996) if we show that the bracketing number $N_{[\,]}(\epsilon,\mathcal F, L_1(P))$ is finite for each $\epsilon>0$. 
To this end 
 let $\epsilon>0$ and for the moment fix $s\in[0,1]$, $\t\in\Theta$ and $z\in C_\tau$. 
Choose $\delta>0$ corresponding to $\epsilon$ as in assumption \ref{F2}.

\noindent Partition $[0,1]$ into finitely many intervals $[s_{j},s_{j+1}]$ such that $F_X(s_{j+1})-F_X(s_{j})\leq \epsilon$ for all $j$. For the fixed $s$, denote the interval containing $s$ by $[s_{j},s_{j+1}]=[s^\ell,s^u]$. 

\noindent  Now choose a finite sup-norm bracketing of length $\gamma$ for the class $\mathcal{L}_S=\{\Lambda_\t|_S:\t\in\Theta\}$  according to (\ref{bracketing}) with $\gamma$ as in assumption \ref{L2} corresponding to the above chosen $\delta$. For the fixed $\t$ this gives a bracket $h^{\ell}\leq h_\t\leq h^u$ of sup-norm length $\gamma$. 

\noindent  Choose a finite sup-norm bracketing of length $\delta$ for the class $\mathcal{L}^1_{\tilde S}=\{\Lambda_0\circ\Lambda_\t^{-1}|_{\tilde S}:\t\in\Theta\}$  according to (\ref{bracketing}). For the fixed $\t$ this gives a bracket $V^{\ell}\leq \Lambda_0\circ\Lambda_\t^{-1}\leq V^u$. 

\noindent Then consider the bounded and increasing function
$$D(z)=\int F_0(V^\ell(z+h^\ell(x))-h_0(x))f_X(x)\,dx$$
and choose a finite partition of the compact $C_\tau$ in intervals $[z_{k},z_{k+1}]$ such that $D(z_{k+1})-D(z_k)<\epsilon$. For the fixed $z$, denote the interval containing $z$ by $[z_{k},z_{k+1}]=[z^\ell,z^u]$.

Now  for the function $f\in \mathcal F$ that is determined by $\t$, $s$ and $z$, a bracket is given by $[f^\ell,f^u]$ with  
\begin{eqnarray*}
f^\ell(x,y) &=& I\{\Lambda_0(y)\leq V^\ell(z^\ell+h^\ell(x))\}I\{x\leq s^\ell\}\\
f^u(x,y) &=& I\{\Lambda_0(y)\leq V^u(z^u+h^u(x))\}I\{x\leq s^u\}
\end{eqnarray*}
with $L_1(P)$-norm
\begin{eqnarray*}
&& \mathbb{E}[I\{\Lambda_0(Y_{i})\leq V^u(z^u+h^u(X_i))\}I\{X_i\leq s^u\}]- \mathbb{E}[I\{\Lambda_0(Y_{i})\leq V^\ell(z^\ell+h^\ell(X_i))\}I\{X_i\leq s^\ell\}]\Big)\\
&\leq& F_X(s^u)-F_X(s^\ell)\\
&&{}+ \int \Big| F_0\left(V^u(z^u+h^u(x))-h_0(x)\right)-F_0\left(V^\ell(z^\ell+h^\ell(x))-h_0(x)\right)\Big|f_X(x)\,dx
\\
 &\leq& 2\epsilon
+\int \Big| F_0\left(V^u(z^u+h^u(x))-h_0(x)\right)
-F_0\left(\Lambda_0(\Lambda_\t^{-1}(z^u+h^u(x)))-h_0(x)\right)\Big|f_X(x)\,dx\\
&&{}+\int \Big| F_0\left(\Lambda_0(\Lambda_\t^{-1}(z^u+h^\ell(x)))-h_0(x)\right)
-F_0\left(V^\ell(z^u+h^\ell(x))-h_0(x)\right)\Big|f_X(x)\,dx\\
&&{}+ \int \Big| F_0\left(\Lambda_0(\Lambda_\t^{-1}(z^u+h^u(x)))-h_0(x)\right)
-F_0\left(\Lambda_0(\Lambda_\t^{-1}(z^u+h^\ell(x)))-h_0(x)\right)\Big|f_X(x)\,dx\\
 &\leq& 4\epsilon
\end{eqnarray*}
by the definition of $[s^\ell,s^u]$ and $[z^\ell,z^u]$ and
using the construction of brackets above (note that $\|V^u-\Lambda_0\circ\Lambda_\t^{-1}\|_\infty\leq \delta$, $\|\Lambda_0\circ\Lambda_\t^{-1}-V^\ell\|_\infty\leq \delta$, $\|h^u-h^\ell\|_\infty\leq\gamma$   and recall assumptions \ref{F2} and \ref{L2}).

There are finitely many such brackets to cover $\mathcal F$ and thus the assertion follows. 
 \hfill $\Box$

\begin{lemma}\label{lem1-random} Under the assumptions of Theorem \ref{theo} (i),
$$\sup_{\t\in\Theta} \|G(\t,h_{\t})-G(\t,\hat h_{\t})\| =o_P(1).$$
\end{lemma}

{\bf Proof.}
According to assumption \ref{N1}  it suffices to show
$$\sup_{\t\in\Theta}\sup_{s\in [0,1]\atop y\in C} |G(\t,h_{\t})(y,s)-G(\t,\hat h_{\t})(y,s)| =o_P(1).$$
Recalling the definition of $G$ in (\ref{G-random}) we see that the assertion follows from Lemma \ref{lemma-smoothing} and uniform continuity of $F_0$ and of $\Lambda_0\circ\Lambda_\t^{-1}$ (uniformly in $\t$). 
 \hfill $\Box$

\section{Identifiability of the model in the random design case}\label{ident-proof}
	
\noindent We prove the assertion of Remark \ref{identifiability}. 
First note that $\varepsilon(\vartheta_1)$ is independent of $X$, and thus  the conditional distribution of $\varepsilon(\vartheta_1)$, i.e.
\begin{eqnarray*}
\mathbb{P}(\varepsilon(\vartheta_1)\leq y\mid X=x)&=& \mathbb{P}(Y\leq \Lambda_{\vartheta_1}^{-1}(y+h_{\vartheta_1}(x)))\mid X=x)\\
&=& F_0(\Lambda_{\vartheta_0} (\Lambda_{\vartheta_1}^{-1}(y+ h_{\vartheta_1}(x)))-h_{\vartheta_0}(x)) 
\end{eqnarray*}
does not depend on $x$. Further, $h_{\vartheta_0}=\Lambda_{\vartheta_0} \circ\Lambda_{\vartheta_1}^{-1}\circ h_{\vartheta_1}$, and for $y\leq 0$ we have $\Lambda_{\vartheta_0} (\Lambda_{\vartheta_1}^{-1}(y+h_{\vartheta_1}(x)))\leq \Lambda_{\vartheta_0} (\Lambda_{\vartheta_1}^{-1}(h_{\vartheta_1}(x)))$ because $\Lambda_{\vartheta_0} \circ\Lambda_{\vartheta_1}^{-1}$ is strictly increasing.
As $F_0$ is strictly increasing by assumption it follows that 
$$H^{-1}(y+H(h_{\vartheta_0}(x)))- h_{\vartheta_0}(x)$$
does not depend on $x$ for $y\in (-\infty,0]$ and $x\in [0,1]$, where for ease of presentation  write $H:= \Lambda_{\vartheta_1}\circ \Lambda_{\vartheta_0}^{-1}$. Thus
$$H^{-1}(y+H(a))-a=H^{-1}(y+H(b))-b$$
for all $y\leq 0$, $a,b\in h_{\vartheta_0}([0,1])$. 
Because $Y$ may take the value $0$ by assumption and $\varepsilon\leq 0$ one obtains $h_{\vartheta_0}([0,1])\cap \mathbb{R}_0^+\neq\emptyset$. 
To conclude the proof we distinguish two cases. 

(1) Let $h_{\vartheta_0}([0,1])\cap \mathbb{R}_0^+=\{0\}$. Set $a=0$, since by assumption $\Lambda_{\vartheta}(0)=0$ for all $\vartheta \in \Theta$, then 
$$H^{-1}(y)=H^{-1}(y+H(b))-b$$
for all $y\leq 0$, $b\in h_{\vartheta_0}([0,1])\subset\mathbb{R}_0^-$. 
Set $c=H^{-1}(y+H(b))$, then it follows that
$$H(c)-H(b)=H(c-b)$$
for all $b,c\in (-\delta,0]$ for some $\delta>0$ and from the assumptions it follows that $\vartheta_1=\vartheta_0$ with $H=\mbox{id}$.

(2) Let $h_{\vartheta_0} ([0,1])\cap \mathbb{R}_0^+=I$ be an interval of positive length. For $a\in I$ one has $y:=-H(a)\leq 0$, and  
$$0=H^{-1}(0)=a+H^{-1}(-H(a)+H(b))-b$$
and thus 
$$H(b-a)=H(b)-H(a)$$
for all $a,b\in I$. 
From the assumptions it follows that $\vartheta_1=\vartheta_0$ with $H=\mbox{id}$ and thus identifiability of the model.

\hfill $\Box$

\section*{References}

\begin{description}

\item Box, G. E. P. and Cox, D. R. (1964).
An analysis of transformations.
{\it J. Roy. Statist. Soc. Ser. B}  {\bf 26}, 211--252.

\item Brown, L.D. and Low M.G. (1996). Asymptotic equivalenfce of nonparametric regression and white noise. \textit{The annals of statistics} \textbf{24}, 2384--2398.

\item  Carroll, R. J. and Ruppert, D. (1988).
\textit{Transformation and weighting in regression}. 
Monographs on Statistics and Applied Probability. 
Chapman \& Hall, New York.


\item Colling, B. and Van Keilegom, I. (2016).
Goodness-of-fit tests in semiparametric transformation models. 
 TEST \textbf{25}, 291--308. 

\item Daouia, A., Noh, H. and Park, B. U. (2016).  Data envelope fitting with constrained polynomial splines.
\textit{J. R. Stat.  Soc.  B.} \textbf{78}, 3--30.
\item Drees, H., Neumeyer, N. and Selk, L. (2018). Estimation and hypotheses tests in boundary regression models. \textit{Bernoulli}, to appear.  \\
http://www.bernoulli-society.org/index.php/publications/bernoulli-journal

\item Girard, S.\ and Jacob, P. (2008). Frontier estimation via kernel regression on high power-transformed data. \textit{J. Multivariate Anal.} 
\textbf{99}, 403--420.

\item Hall, P. and Van Keilegom, I. (2009). Nonparametric "regression" when errors are positioned at end-points. \textit{Bernoulli} \textbf{15}, 614--633.

\item Hall, P., Park, B.U. and  Stern, S.E. (1998). On polynomial estimators of frontiers and boundaries. \textit{J. Multivariate Anal.}  \textbf{66}, 71--98.

\item Hõrdle, W., Park, B.U. and Tsybakov, A.B. (1995). Estimation of a non sharp support boundaries. \textit{J. Multivariate Anal.} \textbf{55}, 205--218.

\item Horowitz, J. L. (2009). \textit{Semiparametric and nonparametric methods in econometrics}. Springer Series in Statistics. Springer, New York.

\item Jirak, M., Meister, A.\ and Rei{\ss}, M. (2013). Asymptotic equivalence for nonparametric regression with non-regular errors. \textit{Probability theory and relative fields} {\bf 155}, 201--229.

\item Jirak, M., Meister, A.\ and Rei{\ss}, M. (2014). Adaptive estimation in nonparametric regression with one-sided errors. \textit{Ann.\ Statist.} {\bf 42}, 1970--2002.

\item Jones, M.C. and Pewsey, A. (2009). Sinh-arcsinh distributions. \textit{Biometrika} \textbf{96}, 761--780.

\item  Kosorok, M. R. (2008). \textit{Introduction to empirical processes and semiparametric inference.} Springer Series in Statistics. Springer, New York.

\item Linton, O., Sperlich, S. and Van Keilegom, I. (2008).
Estimation on a semiparametric transformation model.
{\it Ann. Statist.}  {\bf 36}, 686--718.

\item Manly, B.F.J (1976). Exponential data transformations. \textit{Journal of the Royal Statistical Society. Series D} \textbf{25}, 37--42.

\item Markowitz, H. (1952). Portfolio Selection. \textit{The Journal of Finance} \textbf{ 7}, 77--91.

\item Meister, A.\ and  Rei{\ss}, M.  (2013). Asymptotic
equivalence for nonparametric regression with non-regular errors. \textit{Probab.\ Th.\ Rel.\ Fields} {\bf 155}, 201--229.

\item Mu, Y. and He, X. (2007).
Power transformation toward a linear regression quantile.
{\it J. Amer. Statist. Assoc.}  {\bf 102}, 269--279.

 \item M\"{u}ller, U.U.\ and Wefelmeyer W. (2010). Estimation in nonparametric regression with non-regular errors. \textit{Comm.\ Statist.\ Theory Methods} {\bf 39}, 1619--1629.

\item Pollard, D. (1984). \textit{Convergence of Stochastic Processes}. Springer, New York.

\item Powell, J. (1991). Estimation of monotonic regression models under quantile restrictions. In:
Barnett, W., Powell, J., and Tauchen, G. (eds.), {\em Nonparametric and Semiparametric Methods in
Econometrics}, 357--384. Cambridge University Press, New York.

\item Simar, L. and Wilson, P.W. (1998). Sensitivity analysis of efficiency scores: how to bootstrap in nonparametric frontier models. {\em Management Science} {\bf 44}, 49--61.

 \item van der Vaart, A.W. and Wellner, J.A. (1996). \textit{Weak convergence and empirical processes}. Springer, New York.

\item Wilson, P.W. (2003). Testing independence in models of productive efficiency. {\em Journal of Productivity Analysis} {\bf 20}, 361--390.

\item Yeo, I-K. and Johnson, R. A. (2000).
A new family of power transformations to improve normality or symmetry. 
{\it Biometrika}  {\bf 87}, 954--959.
\end{description}

\end{appendix}

\newpage

\begin{minipage}{\textwidth}
\center\epsfig{file=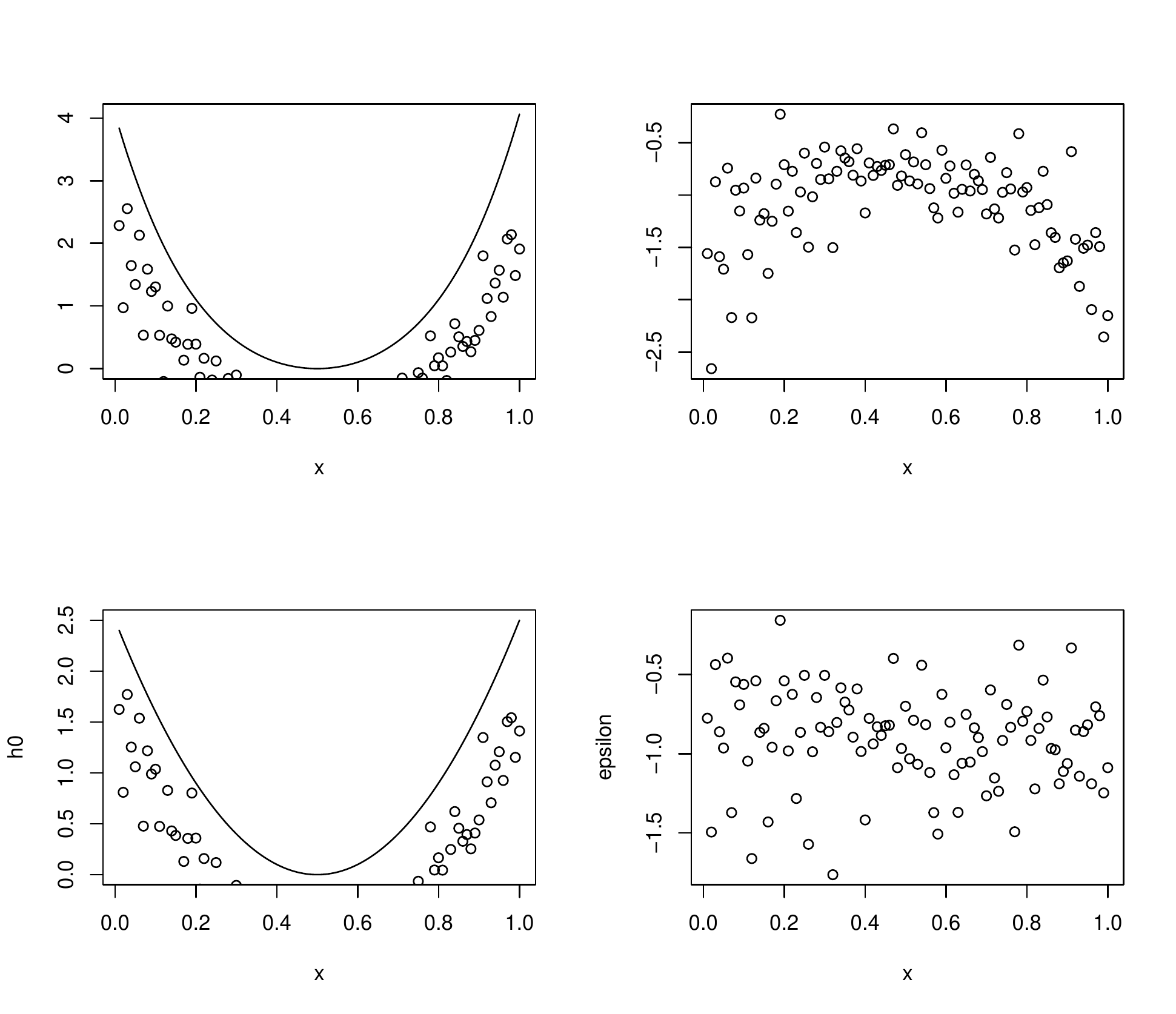, height=0.4\textheight}
\captionof{figure}{Original data (upper panel) and transformed data (\ref{model0}) (lower panel) with $h_0(x)= 10(x-\textstyle{\frac 12})^2$, $n=100$ and $\varepsilon_{i,n}\sim \mbox{Weibull}(1,3)$ with Yeo and Johnson transformation $\Lambda_{0.5}$ as defined in Example \ref{yeo-johnson}. The design points are equidistant. The left panels show the data and regression functions, the right panels show the errors. 
}\label{graphic1}

\center\epsfig{file=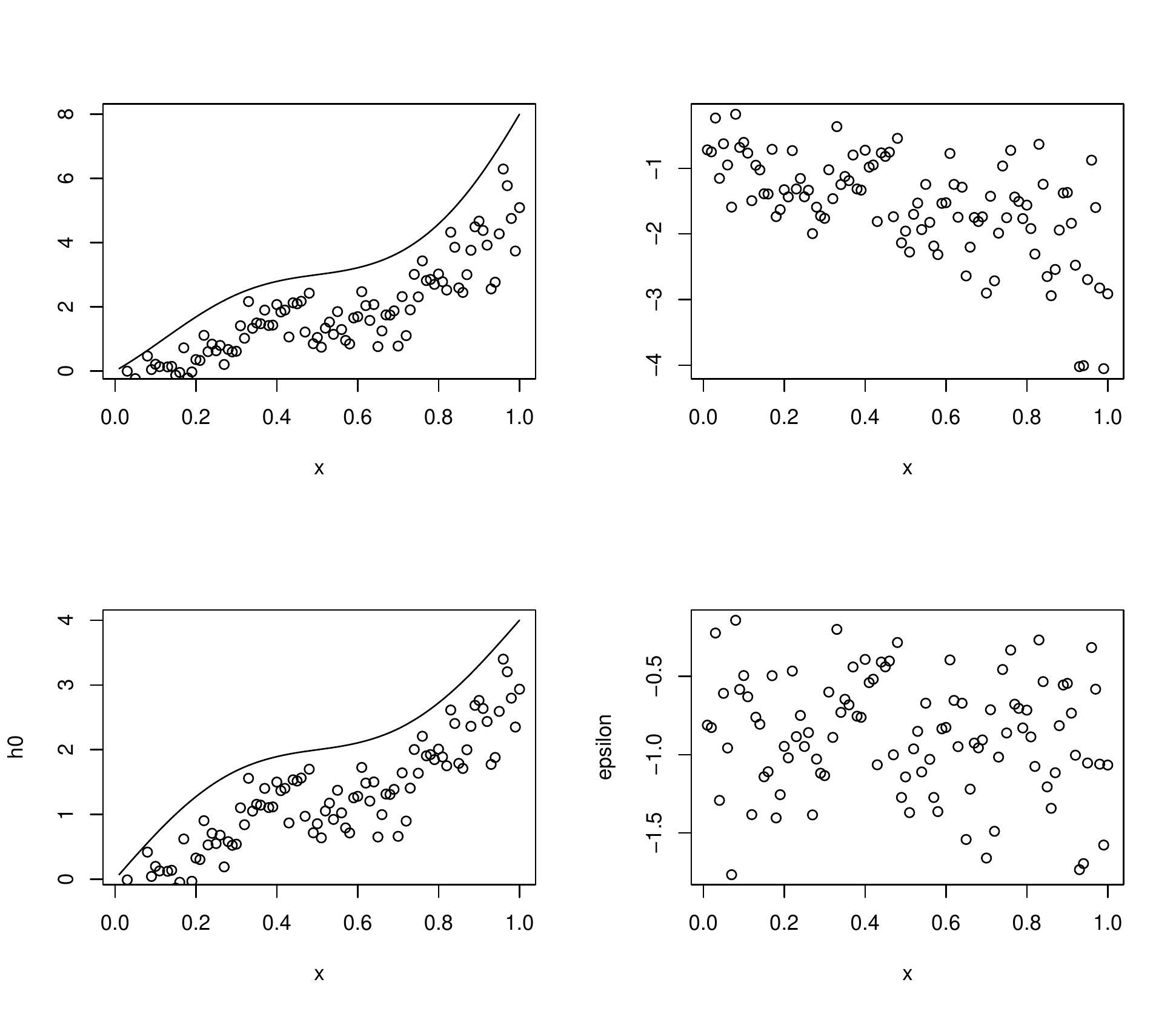,height=0.4\textheight}
\captionof{figure}{The setting is similar to Figure \ref{graphic1} but with $h_0(x)= \textstyle{\frac 12}\sin(2\pi x)+4x$.
}\label{graphic2}
\end{minipage}

\begin{minipage}[t]{0.48\textwidth}
	\center\epsfig{file=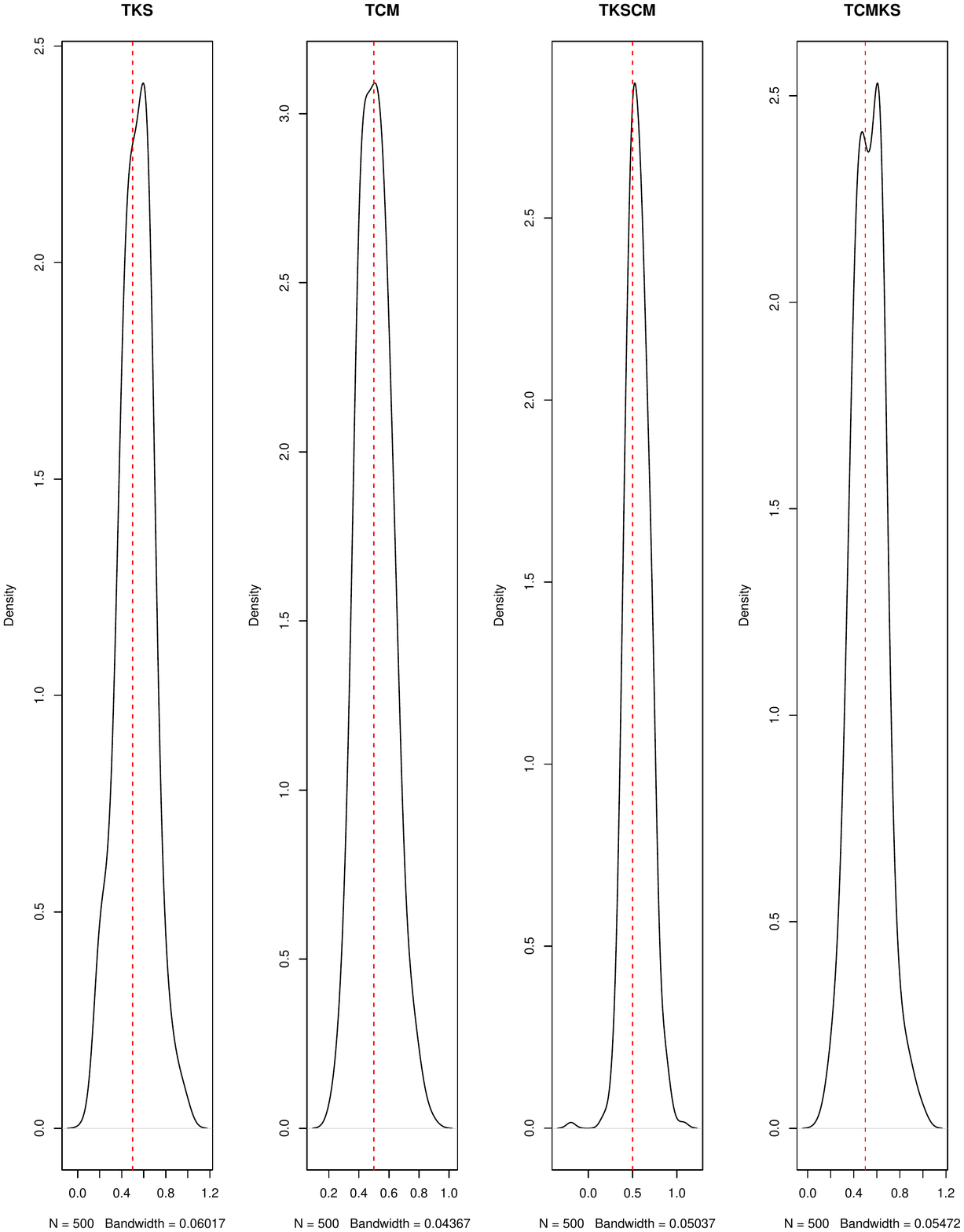,width=\textwidth}
	\captionof{figure}{Density function of the four estimators TKS, TCM, TCMKS and TKSCM for the Model (\ref{simu1}) with a sample size $n=100$ and  bandwidths $a_n=b_n/20$ with $b_n=n^{-1/3}$. This corresponds to results in Table \ref{Table2}. The vertical dashed line corresponds the true parameter $\t_0=0.5$.}\label{est}
\end{minipage}
\ \ \ 
\begin{minipage}[t]{0.48\textwidth}
	\center\epsfig{file=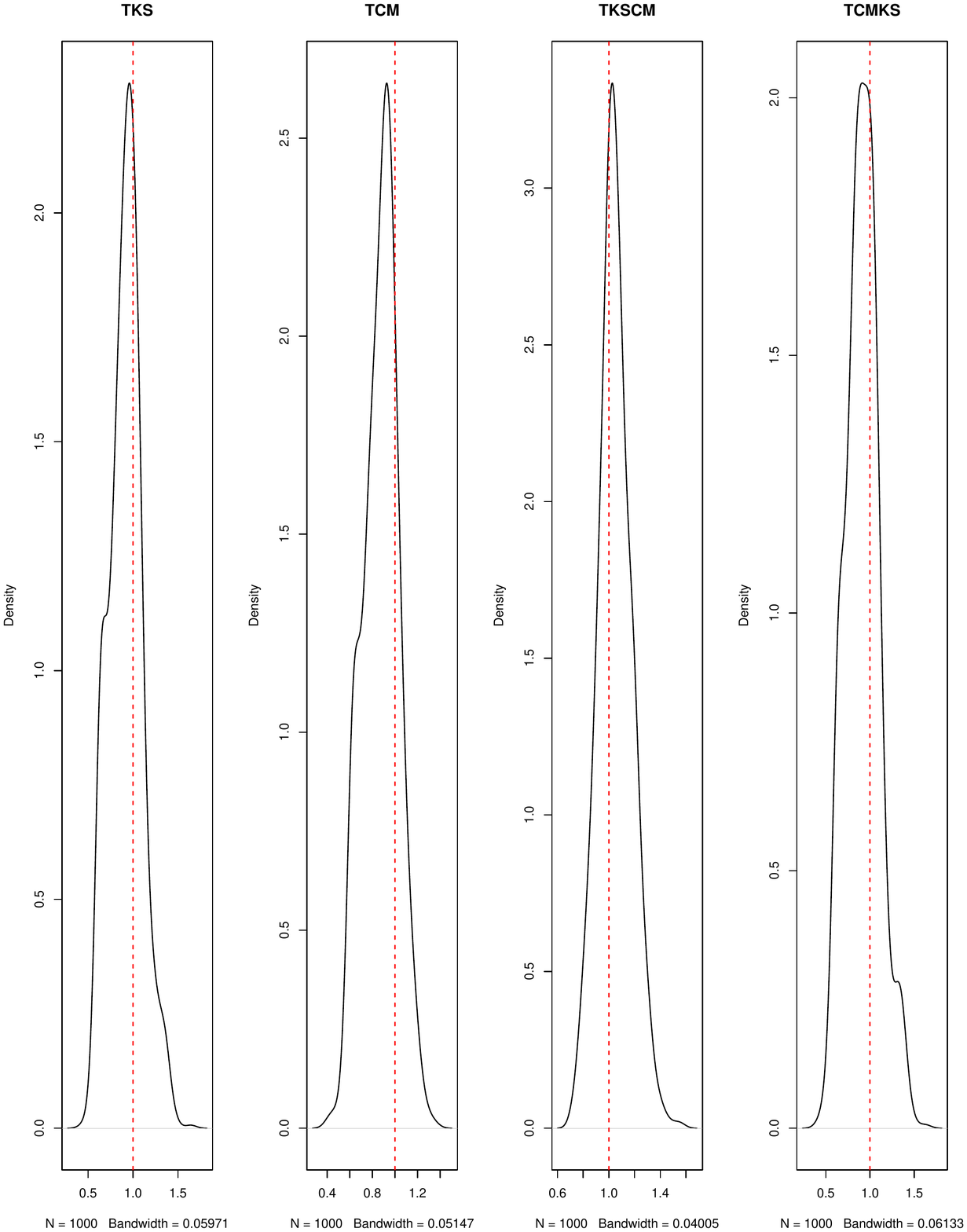,width=\textwidth}
	\captionof{figure}{Density function of the four estimators TKS, TCM, TCMKS and TKSCM for the Model (\ref{simu4}) with a sample size $n=100$ and  bandwidths $a_n=b_n/2$ with $b_n=n^{-1/3}$. This corresponds to results in Table \ref{Table7}. The vertical dashed line corresponds the true parameter $\t_0=1$.}\label{est2}
\end{minipage}

\begin{table}	
		\begin{tabular}{|l|c|c|c|c|}
			\hline
			$n=50$& TKS & TCM & TKSCM & TCMKS \\
			\hline
			$\t_0=0$
			& 0.162 0.162 (0.120)  & \textbf{0.095} 0.122 (\underline{0.060}) & 0.216 0.241 (0.149) & 0.196 0.183 (0.141)  
			\\
			\hline
			$\t_0=0.5$
			& 0.643 0.646 (0.142) & \textbf{0.576} 0.595 (\underline{0.080}) & 0.766 0.800 (0.250) & 0.691 0.646 (0.204)
			\\
			\hline
			$\t_0=1$
			& 1.120 1.190 (0.232) & \textbf{1.100} 1.090 (\underline{0.121}) & 1.340 1.380 (0.335) & 1.200 1.310 (0.287)
			\\
			\hline
			$\t_0=1.5$
			& \textbf{1.610} 1.720 (0.228) & 1.640 1.670 (\underline{0.133}) & 1.900 1.920 (0.336) & 1.620 1.780 (0.282)
			\\
			\hline
			$\t_0=2$
			& 1.810 2.020 (0.357) & \textbf{2.110} 2.130 (\underline{0.076}) & 2.310 2.460 (0.179) & 1.860 2.060 (0.297)
			\\
			\hline
		\end{tabular}
		
		\bigskip
		
		\begin{tabular}{|l|c|c|c|c|}
			\hline
			$n=100$ & TKS & TCM & TKSCM & TCMKS \\
			\hline
			$\t_0=0$
			& 0.014 -0.039 (0.055) & -0.014 -0.029  (\underline{0.019}) & 0.098 0.092 (0.031) & \textbf{-0.006} -0.021 (0.026)
			\\
			\hline
			$\t_0=0.5$
			& 0.483 0.496 (0.041) & 0.461 0.451 (\underline{0.026}) & 0.625 0.618 (0.049) & \textbf{0.503} 0.523 (0.047)
			\\
			\hline
			$\t_0=1$
			& 0.951 0.964 (0.062) & 0.949 0.950 (\underline{0.034}) & 1.150 1.140 (0.059) & \textbf{1.000} 1.010 (0.071)
			\\
			\hline
			$\t_0=1.5$
			& \textbf{1.500} 1.490 (0.050) & 1.470 1.450 (\underline{0.040}) & 1.670 1.660 (0.078) & 1.520 1.520 (0.077)
			\\
			\hline
			$\t_0=2$
			& 1.960 2.000 (0.071) & \textbf{1.970} 1.960 (\underline{0.034}) & 2.150 2.150 (0.065) & 1.970 2.030 (0.074)
			\\
			\hline
		\end{tabular}
		\captionof{table}  {
		\it Mean, median and MISE for Model (\ref{simu1}) for $n=50$ and $n=100$ with $a_n=b_n/2$. }\label{Table1}
\end{table}

\begin{table}
		\begin{tabular}{|l|c|c|c|c|}
			\hline
			$n=50$& TKS & TCM & TKSCM & TCMKS \\
			\hline
			$\t_0=0$
			& 0.154 0.146 (0.126)  & \textbf{-0.059} -0.051 (\underline{0.050}) & 0.198 0.205 (0.114) & 0.173 0.171 (0.128)  
			\\
			\hline
			$\t_0=0.5$
			& 0.613 0.635 (0.118) & \textbf{0.511} 0.513 (\underline{0.066}) & 0.717 0.735 (0.173) & 0.658 0.646 (0.163)
			\\
			\hline
			$\t_0=1$
			& 1.120 1.190 (0.182) & \textbf{1.030} 1.030 (\underline{0.091}) & 1.280 1.300 (0.230) & 1.190 1.270 (0.221)
			\\
			\hline
			$\t_0=1.5$
			& 1.630 1.690 (0.164) & \textbf{1.560} 1.540 (\underline{0.095}) & 1.830 1.810 (0.233) & 1.620 1.740 (0.209)
			\\
			\hline
			$\t_0=2$
			& 1.880 2.040 (0.272) & \textbf{2.050} 2.080 (\underline{0.066}) & 2.280 2.370 (0.149) & 1.870 2.040 (0.272)
			\\
			\hline
		\end{tabular}
		
		\bigskip
		
		\begin{tabular}{|l|c|c|c|c|}
			\hline
			$n=100$ & TKS & TCM & TKSCM & TCMKS \\
			\hline
			$\t_0=0$
			& \textbf{0.003} 0.046 (0.040) &  0.045 0.057 (\underline{0.020}) & 0.078 0.067 (0.026) & 0.005 0.022 (0.025)
			\\
			\hline
			$\t_0=0.5$
			& 0.454 0.448 (0.037) & 0.418 0.406 (\underline{0.028}) & 0.581 0.566 (0.033) & \textbf{0.473} 0.478 (0.046)
			\\
			\hline
			$\t_0=1$
			& 0.938 0.950 (0.053) & 0.913 0.918 (\underline{0.038}) & 1.100 1.080 (0.047) & \textbf{0.984} 0.986 (0.070)
			\\
			\hline
			$\t_0=1.5$
			& 1.450 1.420 (0.052) & 1.410 1.390 (\underline{0.041}) & 1.600 1.590 (0.050) & \textbf{1.470} 1.440 (0.062)
			\\
			\hline
			$\t_0=2$
			& \textbf{1.950} 1.960 (0.057) & 1.930 1.910 (\underline{0.037}) & 2.100 2.100 (0.053) & 1.940 1.980 (0.066)
			\\
			\hline
		\end{tabular}
	\captionof{table}  {
		\it Mean, median and MISE for Model (\ref{simu1}) for $n=50$ and $n=100$ with $a_n=b_n/20$. }\label{Table2}
\end{table}

\begin{table}
	\begin{center}
		\begin{tabular}{|l|c|c|c|c|}
			\hline
			$n=50$& TKS & TCM & TKSCM & TCMKS \\
			\hline
			$\t_0=0$
			& 0.054 0.062 (0.218) & -0.029 -0.029 (0.005) & \textbf{0.000} 0.000 (\underline{0.000}) & 0.050 0.046 (0.182) 
			\\
			\hline
			$\t_0=0.5$
			& 0.419 0.427 (0.038) & 0.412 0.408 (0.027) & \textbf{0.469} 0.478 (\underline{0.025}) & 0.458 0.470 (0.030)
			\\
			\hline
			$\t_0=1$
			& 0.887 0.881 (0.080) & 0.791 0.770 (0.104) & \textbf{0.942} 0.971 (0.077) & 0.934 0.933 (\underline{0.069}) 
			\\
			\hline
			$\t_0=1.5$
			& 1.360 1.350 (0.130) & 1.210 1.270 (0.214) & 1.310 1.420 (0.319) & \textbf{1.410} 1.350 (\underline{0.111})
			\\
			\hline
			$\t_0=2$
			& 1.800 1.790 (0.165) & 1.560 1.690 (0.504) & 1.249 1.769 (1.739) & \textbf{1.810} 1.790 (\underline{0.141})
			\\
			\hline
		\end{tabular}
		
		\bigskip
		
		\begin{tabular}{|l|c|c|c|c|}
			\hline
			$n=100$ & TKS & TCM & TKSCM & TCMKS \\
			\hline
			$\t_0=0$
			& 0.025 0.043 (0.004) & -0.037 -0.038 (\underline{0.003}) & \textbf{0.012} 0.009 (\underline{0.002}) & -0.034 -0.054 (0.004)
			\\
			\hline
			$\t_0=0.5$
			& 0.459 0.458 (0.017) & 0.406 0.405 (0.016) & \textbf{0.520} 0.518 (\underline{0.008}) & 0.449 0.453 (0.017)
			\\
			\hline
			$\t_0=1$
			& 0.922 0.945 (0.038) & 0.845 0.856 (0.046) & \textbf{1.030} 1.040 (\underline{0.017}) & 0.907 0.909 (0.044)
			\\
			\hline
			$\t_0=1.5$
			& 1.430 1.350 (0.056) & 1.290 1.310 (0.080) & \textbf{1.540} 1.540 (\underline{0.032}) & 1.400 1.350 (0.063)
			\\
			\hline
			$\t_0=2$
			& \textbf{1.930} 1.970 (\underline{0.072}) & 1.740 1.770 (0.126) & 1.880 2.010 (0.419) & 1.850 1.790 (0.090)
			\\
			\hline
		\end{tabular}
	\captionof{table}  {
		\it Mean, median and MISE for Model (\ref{simu4}) for $n=50$ and $n=100$ with $a_n=b_n/2$. }\label{Table7}
	\end{center}
\end{table}

\begin{table}
	\begin{center}
		\begin{tabular}{|l|c|c|c|c|}
			\hline
			$n=50$& TKS & TCM & TKSCM & TCMKS \\
			\hline
			$\t_0=0$
			& 0.021 0.062 (0.177) & 0.053 0.057 (\underline{0.007}) & \textbf{0.013} 0.010 (\underline{0.007}) & 0.017 0.062 (0.141) 
			\\
			\hline
			$\t_0=0.5$
			& 0.404 0.415 (0.044) & 0.375 0.380 (0.034) & \textbf{0.449} 0.461 (\underline{0.029}) & 0.434 0.445 (0.034)
			\\
			\hline
			$\t_0=1$
			& 0.851 0.830 (0.086) & 0.737 0.695 (0.119) & \textbf{0.915} 0.931 (\underline{0.074}) & 0.895 0.891 (0.076) 
			\\
			\hline
			$\t_0=1.5$
			& 1.330 1.350 (0.150) & 1.130 1.200 (0.260) & 1.310 1.410 (0.293) & \textbf{1.380} 1.350 (\underline{0.119})
			\\
			\hline
			$\t_0=2$
			& 1.760 1.790 (0.176) & 1.520 1.610 (0.435) & 1.409 1.770 (1.250) & \textbf{1.770} 1.790 (\underline{0.156})
			\\
			\hline
		\end{tabular}
		
		\bigskip
		
		\begin{tabular}{|l|c|c|c|c|}
			\hline
			$n=100$ & TKS & TCM & TKSCM & TCMKS \\
			\hline
			$\t_0=0$
			& 0.027 0.043 (0.004) & -0.038 -0.040 (\underline{0.003}) & \textbf{0.011} 0.008 (\underline{0.002}) & 0.037 0.062 (0.004)
			\\
			\hline
			$\t_0=0.5$
			& 0.462 0.466 (0.018) & 0.411 0.407 (0.015) & \textbf{0.524} 0.522 (\underline{0.008}) & 0.453 0.460 (0.019)
			\\
			\hline
			$\t_0=1$
			& 0.930 0.950 (0.037) & 0.849 0.856 (0.043) & \textbf{1.040} 1.030 (\underline{0.018}) & 0.917 0.915 (0.042)
			\\
			\hline
			$\t_0=1.5$
			& 1.440 1.390 (0.058) & 1.300 1.320 (0.078) &  \textbf{1.550} 1.540 (\underline{0.031}) & 1.420 1.350 (0.065)
			\\
			\hline
			$\t_0=2$
			&  \textbf{1.930} 1.960 (\underline{0.067}) & 1.730 1.760 (0.130) & 1.830 2.010 (0.509) & 1.840 1.790 (0.091)
			\\
			\hline
		\end{tabular}
		\captionof{table}  {
		\it Mean, median and MISE for Model (\ref{simu4}) for $n=50$ and $n=100$ with $a_n=b_n/20$. }\label{Table8}
\end{center}
\end{table}

\begin{table}
	\begin{center}
		\begin{tabular}{|l|c|c|c|}
			\hline
			Method & Pearson & Kendall & Spearman \\
			\hline
			Original data
			& -0.273 & -0.165 & -0.234   
			\\
			\hline
			True parameter $\vartheta_0$
			& 0.005 & 0.003 & 0.004
			\\
			\hline
			TKS
			& 0.008 & 0.004 & 0.007
			\\
			\hline
			TCM
			& 0.024 & 0.014 & 0.021
			\\
			\hline
			TKSCM
			& 0.011 & 0.007 & 0.009
			\\
			\hline
			TCMKS & 0.003 & 0.001 & 0.002
			\\
			\hline
		\end{tabular}
		\captionof{table}  {
		\it Pearson's, Kendall's and Spearman's correlation coefficients (the average over $1000$ iterations) between the covariates and the errors for the model \eqref{simu4} when $n=100$. The first line corresponds to the correlations for the original data while the second line is for the true transformation parameter ($\vartheta_0=0.5$). The last four lines correspond to the correlations for each estimator.}\label{TableCor2paper}
\end{center}
\end{table}

 \newpage
 ~
 \newpage

\begin{center}

\begin{Large}
 \textbf{ Supplementary material to ``Semi-parametric transformation boundary regression models"
 }
 
\end{Large}

\end{center}

\bigskip

\begin{appendix}

\setcounter{section}{2}

\section{Proofs of asymptotic results in the fixed design case}

\subsection{Proof of Lemma \ref{consistency}}\label{proof-consistency}
	
\noindent To prove Lemma \ref{consistency}, we first need the following technical lemma.

\begin{lemma}\label{lem:unif_consist}
	Assume model (\ref{model0}) holds under assumptions \ref{F1'}, \ref{X1'} and \ref{B1'}. Then we have
	\begin{eqnarray*}\label{Lemma2dot2}
	\sup_{x \in [0,1]}  \min_{ \substack{ i \in \{1,\ldots,n\} \\ |x_{i,n}-x|\leq b_n }} |\varepsilon_{i,n}|=o_P(1).
	\end{eqnarray*}
\end{lemma}

{\bf Proof.} The proof is similar to the proof of Lemma A.2 in Drees et al. (2018) but some adaptations are needed to deal with non-equidistant fixed design points. Let $Z_1,Z_2,\dots$ be iid with the same distribution as $-\varepsilon_{i,n}$ with cumulative distribution function $U$.
To prove the result, we shall show that 
	\begin{equation*}
	\lim_{n \to \infty} \mathbb{P}\left( \sup_{x \in [0,1]}  \min_{ \substack{ i \in \{1,\ldots,n\} \\ |x_{i,n}-x|\leq b_n }} Z_i> \epsilon \right)=0, \ \ \ \epsilon>0.
	\end{equation*}	
	For $n \geq 1$, let $0<k \leq n$, $x \in [0,1]$ and set $I_n=[x-b_n,x+b_n]$. Assume that exactly $k$ points lie in $I_n$, say
	\begin{eqnarray*}
	x_{m+1,n} < \dots <x_{m+k,n} \in I_n 
	\end{eqnarray*}
	for some $m < n+1-k$. 
	We shall distinguish two cases.
	\begin{enumerate}[label=(\textbf{\arabic{*}})]
		\item \label{1} If  $(x_{m,n},x_{m+k+1,n}) \in [0,1]^2$, it means that 
		\begin{eqnarray*}
		2b_n=|I_n| < x_{m+k+1,n}-x_{m,n}
		= \sum_{j=m}^{m+k} \left(x_{j+1,n}-x_{j,n}\right)
		\leq (k+1)\bar{\Delta}_n,
		\end{eqnarray*}
		since $\bar{\Delta}_n \geq x_{j,n}-x_{j-1,n}$ for any $1 \leq j \leq n+1$.
		\item \label{2} If $x_{m,n}$ or $x_{m+k+1,n} $ do not exist, which means that either $x_{m+1,n}=x_{0,n}=0$ or $x_{m+k+1,n}=x_{n+1,n}=1$. Consider the first case $x_{m+1,n}=x_{0,n}$ (the extremal case is $x=0$). Then we have  
		\begin{eqnarray*}
		b_n=\frac{|I_n|}{2} < x_{k,n}-x_{0,n}
		= \sum_{j=0}^{k-1} \left(x_{j+1,n}-x_{j,n}\right)
		\leq k\bar{\Delta}_n.
		\end{eqnarray*}
		A similar inequality holds for $x_{m+k+1,n}=x_{n+1,n}=1$ (with the extremal case $x=1$).	
	\end{enumerate}
	In both cases, \ref{1} and \ref{2} yield to 
	\begin{eqnarray*}
	b_n< k\bar{\Delta}_n \Rightarrow k> \frac{b_n}{\bar{\Delta}_n}, \ \ \ n \geq 1.
	\end{eqnarray*}
	
Then, for all $y>0$, we have with $d_n:=\lceil \frac{b_n}{ \bar \Delta_n} \rceil $
\begin{eqnarray*}
	\mathbb{P}\left( \sup_{x \in [0,1]}  \min_{ \substack{ i \in \{1,\ldots,n\} \\ |x_{i,n}-x|\leq b_n }} Z_{i}> y \right  )  &\leq&   \mathbb{P}\left( \left\lbrace  \max_{j \in \{1,\ldots,n-d_n\}} \min_{i \in \{j,\ldots,j+d_n\}} Z_{i}> y \right\rbrace	\right) \\
	&\leq& \sum_{j=1}^{n-d_n}\mathbb{P}\left(\min_{i\in\{j,\ldots,j+d_n\}}Z_i>y\right)\\
	&=& (n-d_n)\mathbb{P}\left(\min_{i\in\{1,\ldots,d_n+1\}}Z_i>y\right)\\
	&=&(n-d_n)\overline{U}(y)^{d_n+1}.
	\end{eqnarray*}
Thus it remains to show that for all $\epsilon>0$
\[(n-d_n)\overline{U}(\epsilon)^{d_n+1}\nto 0\]
which is true since  $d_n \underset{n \to \infty}{\sim} \frac{b_n}{ \bar \Delta_n}$ and
\beq
\frac{b_n}{ \bar \Delta_n} \log(\overline{U}(\epsilon))+\log(n-d_n)&\leq& \frac{b_n}{ \bar \Delta_n}\log(\overline{U}(\epsilon))+\log(n)\\
&=&\log(n)\left(\frac{b_n}{ \bar \Delta_n \log(n)}\log(\overline{U}(\epsilon))+1\right)\\
&\nto& -\infty
\eeq
since $\overline{U}(\epsilon)<1$ under \ref{F1'} and $\frac{b_n}{\bar \Delta_n \log(n)} \nto \infty$ under \ref{B1'}. This concludes the proof.

 \hfill $\Box$

\bigskip

\noindent
The
{\bf proof of Lemma \ref{consistency}} is analogous to  the proof of Lemma \ref{consistency-random}.

\subsection{Proof of Theorem \ref{theo} in the fixed design case}\label{proof-theo-fixed}

The first part of the proof is similar to the random design case. Here, we use
\begin{eqnarray*}
\sup_{\t\in\Theta}|M_n( \t)-M(\t)|&\leq &  \sup_{\t\in\Theta}\|G_n(\t,\hat h_{ \t})-\bar G_n(\t,\hat h_{\t})\|+
\sup_{\t\in\Theta}\|\bar G_n(\t,\hat h_{ \t})-\tilde G_n(\t,\hat h_{\t})\| \\
&&{}+ \sup_{\t\in\Theta}\|\tilde G_n(\t,\hat h_{\t})-G(\t,\hat h_{\t})\|+
\sup_{\t\in\Theta} \|G(\t,\hat h_{\t})-G(\t, h_{\t})\|, 
\end{eqnarray*}
where the definition for $M$ and $G$ is as in the random case, and
\begin{eqnarray*}
\bar G_n(\t, h)(y,s)&=& \frac 1n\sum_{i=1}^n I\{\Lt(Y_{i,n})-h(\tin)\leq y\} \big( I\{\tin\leq s\}-F_{X}(s)\big).\\
\end{eqnarray*}
Further, 
\begin{eqnarray}
\label{G_n-tilde}
\tilde G_n(\t,h)(y,s)
&=&\frac 1n\sum_{i=1}^n F_0\left(\Lambda_0(\Lt^{-1}(y+h(\tin)))-h_0(\tin)\right) I\{\tin\leq s\}\\
&&{}-F_X(s)\sum_{i=1}^n F_0\left(\Lambda_0(\Lt^{-1}(y+h(\tin)))-h_0(\tin)\right)
\nonumber
\end{eqnarray}
is a Riemann-sum approximation of $G(\t,h)(y,s)$. Note that for any deterministic function $h$ we have $\tilde G_n(\t,h)=\mathbb{E}[\bar G_n(\t,h)]$.
The assertion of the theorem follows from 
\begin{eqnarray}\label{glivenko-cantelli_2}
\sup_{\t\in\Theta}\|G_n(\t,\hat h_{ \t})-\bar G_n(\t,\hat h_{\t})\|\leq \sup_{s\in[0,1]}|\hat F_{X,n}(s)-F_X(s)|=o(1)
\end{eqnarray}
and from 
Lemmas \ref{lem2}--\ref{lem3} by an application of the arg-max theorem.
For (\ref{glivenko-cantelli_2}) note that with assumption \ref{X1''}
\begin{eqnarray}\label{glivenko-cantelli}
\sup_{s\in[0,1]}|\hat F_{X,n}(s)-F_X(s)|
&=&\sup_{s\in[0,1]}\bigg|\frac 1n\sum_{i=1}^n I\{x_{i,n}\leq s\}-\int_0^s f_X(x)\,dx\bigg|\\
&\leq &\sup_{s\in[0,1]}\bigg|\sum_{i=1}^n \int_{x_{i-1,n}}^{x_{i,n}}f_X(x)\,dxI\{x_{i,n}\leq s\}-\int_0^s f_X(x)\,dx\bigg| +o(1)\nonumber\\
&=&\sup_{s\in[0,1]}\bigg|\int_{\max\{x_{i,n}| x_{i,n}\leq s\}}^sf_X(x)\,dx\bigg|+o(1)\nonumber\\
&=&\bar\Delta_n\sup_{x\in[0,1]}f_X(x)+o(1) =o(1).
\end{eqnarray}
 \hfill $\Box$

\begin{lemma}\label{lem2} Under the assumptions of Theorem \ref{theo} (ii),
$$\sup_{\t\in\Theta}\|\bar G_n(\t,\hat h_{ \t})-\tilde G_n(\t,\hat h_{ \t})\| =o_P(1).$$
\end{lemma}

{\bf Proof.} As in the proof of Lemma \ref{lem2-random} we assume in what follows that (\ref{E_n}) holds. 
We only consider the difference between the first sum in the definitions of $G_n(\t, h)$ and the first sum in  $\tilde G_n(\t, h)$ (see (\ref{G_n}) and (\ref{G_n-tilde}), respectively). The difference of the second sums can be treated similarly. 
Applying (\ref{E_n}) the first sum in $ G_n(\t, \hat h_\t)(y,s)$ can be nested as 
\begin{eqnarray*}
&&\frac 1n\sum_{i=1}^n I\{\Lt(Y_{i,n})-h_\t({\tin})\leq y-a_n\}I\{x_{i,n}\leq s\}\\
&\leq& \frac 1n\sum_{i=1}^n I\{\Lt(Y_{i,n})-\hat h_\t({\tin})\leq y\}I\{x_{i,n}\leq s\}\\
&\leq& \frac 1n\sum_{i=1}^n I\{\Lt(Y_{i,n})-h_\t({\tin})\leq y+a_n\}I\{x_{i,n}\leq s\}
\end{eqnarray*}
while the first sum in $ \tilde G_n(\t, \hat h_\t)(y,s)$ can be nested as 
\begin{eqnarray*}
&&\frac 1n\sum_{i=1}^n F_0\left(\Lambda_0(\Lt^{-1}(y-a_n+h_\t(\tin)))-h_0(\tin)\right)I\{x_{i,n}\leq s\}\\
&\leq& \frac 1n\sum_{i=1}^n F_0\left(\Lambda_0(\Lt^{-1}(y+\hat h_\t(\tin)))-h_0(\tin)\right)I\{x_{i,n}\leq s\}\\
&\leq& \frac 1n\sum_{i=1}^n F_0\left(\Lambda_0(\Lt^{-1}(y+a_n+h_\t(\tin)))-h_0(\tin)\right)I\{x_{i,n}\leq s\}.
\end{eqnarray*}
Thus we have to consider
\begin{eqnarray*}
H_{n,\t}^{(1)}(y,s) 
&=&\frac 1n\sum_{i=1}^n \Big(I\{\Lt(Y_{i,n})-h_\t({\tin})\leq y+a_n\} \\
&&\qquad{}-  F_0\left(\Lambda_0(\Lt^{-1}(y+a_n+h_\t(\tin)))-h_0(\tin)\right)\Big)I\{x_{i,n}\leq s\}\\
H_{n,\t}^{(2)}(y,s) 
&=&\frac 1n\sum_{i=1}^n \Big(F_0\left(\Lambda_0(\Lt^{-1}(y+a_n+h_\t(\tin)))-h_0(\tin)\right)\\
&&\qquad{}-F_0\left(\Lambda_0(\Lt^{-1}(y+h_\t(\tin)))-h_0(\tin)\right)\Big)I\{x_{i,n}\leq s\}
\end{eqnarray*}
and the same terms with $y+a_n$ replaced by $y-a_n$, which can be treated completely analogously. We have to show that
$\sup_{\t\in\Theta}\|H_{n,\t}^{(1)}\|=o_P(1)$ and $\sup_{\t\in\Theta}\|H_{n,\t}^{(2)}\|=o(1)$. 




Recall condition \ref{N1} and note that $\sup_{\t\in\Theta}\sup_{s\in [0,1]\atop y\in C}|H_{n,\t}^{(2)}(y,s)|=o(1)$ follows from uniform continuity of $F_0$ and of $\Lambda_0\circ \Lt^{-1}$ uniformly in $\t$  (see \ref{F2} and \ref{L2}), from the representation $h_\t=\Lambda_\t\circ\Lambda_0^{-1}\circ h_0$ and uniform continuity of $\Lambda_\t$ uniformly in $\t$ (see \ref{L1}), and $a_n\to 0$. 

Now to prove $\sup_{\t\in\Theta}\|H_{n,\t}^{(1)}\|=o_P(1)$, let $\epsilon>0$ and for the moment fix $s\in[0,1]$, $\t\in\Theta$ and $y\in C$. 
Choose $\delta>0$ corresponding to $\epsilon$ as in assumption \ref{F2}.
Let $n$ be large enough such that $|a_n|\leq \tau$ for $\tau$ both from \ref{F2} and \ref{L2}. 

\noindent Partition $[0,1]$ into finitely many intervals $[s_{j},s_{j+1}]$ such that $F_X(s_{j+1})-F_X(s_j)<\epsilon$ for all $j$. For the fixed $s$, denote the interval containing $s$ by $[s_{j},s_{j+1}]=[s^\ell,s^u]$. 

\noindent  Now choose a finite sup-norm bracketing of length $\gamma$ for the class $\mathcal{L}_S=\{\Lambda_\t|_S:\t\in\Theta\}$  according to (\ref{bracketing}) with $\gamma$ as in assumption \ref{L2} corresponding to the above chosen $\delta$. For the fixed $\t$ this gives a bracket $h^{\ell}\leq h_\t\leq h^u$ of sunorm length $\gamma$. 

\noindent  Choose a finite sup-norm bracketing  of length $\delta$  for the class $\mathcal{L}^1_{\tilde S}=\{\Lambda_0\circ\Lambda_\t^{-1}|_{\tilde S}:\t\in\Theta\}$ according to (\ref{bracketing}). For the fixed $\t$ this gives a bracket $V^{\ell}\leq \Lambda_0\circ\Lambda_\t^{-1}\leq V^u$. 

\noindent Then consider the bounded and increasing function
$$D_n(y)=\frac{1}{n}\sum_{i=1}^n F_0(V^\ell(y+a_n+h^\ell(\tin))-h_0(\tin))$$
and choose a finite partition of the compact $C$ in intervals $[y_{k},y_{k+1}]$ such that $D_n(y_{k+1})-D_n(y_k)<\epsilon$. For the fixed $y$, denote the interval containing $y$ by $[y_{k},y_{k+1}]=[y^\ell,y^u]$.
Note that the brackets depend on $n$. This is suppressed in the  notation because it is not relevant for the remainder of the proof because the number of brackets is $O(\epsilon^{-1})$, uniformly in $n$. 

Now we can nest as follows 
\begin{eqnarray*}
&&I\{\Lambda_0(Y_{i,n})\leq V^\ell(y^\ell+a_n+h^\ell({\tin}))\}I\{\tin\leq s^\ell\}\\
&\leq&
I\{\Lt(Y_{i,n})-h_\t({\tin})\leq y+a_n\}I\{\tin\leq s\}\\
&=& I\{Y_{i,n}\leq \Lambda_\t^{-1}(y+a_n+h_\t({\tin}))\}I\{\tin\leq s\}\\
&\leq& I\{\Lambda_0(Y_{i,n})\leq V^u(y^u+a_n+h^u({\tin}))\}I\{\tin\leq s^u\},
\end{eqnarray*}
and have
\begin{eqnarray*}
&& \frac 1n\sum_{i=1}^n \Big(\mathbb{E}[I\{\Lambda_0(Y_{i,n})\leq V^u(y^u+a_n+h^u({\tin}))\}I\{\tin\leq s^u\}]\\
&&\qquad\qquad{}- \mathbb{E}[I\{\Lambda_0(Y_{i,n})\leq V^\ell(y^\ell+a_n+h^\ell({\tin}))\}I\{\tin\leq s^\ell\}]\Big)\\
&\leq& \hat F_{X,n}(s^u)-\hat F_{X,n}(s^\ell)\\
&&{}+ \frac 1n\sum_{i=1}^n\Big| F_0\left(V^u(y^u+a_n+h^u(\tin))-h_0(\tin)\right)-F_0\left(V^\ell(y^\ell+a_n+h^\ell(\tin))-h_0(\tin)\right)
\Big|
\\
 &\leq& 2\epsilon+o(1)\\
 &&{}+ \frac 1n\sum_{i=1}^n\Big| F_0\left(V^u(y^u+a_n+h^u(\tin))-h_0(\tin)\right)-F_0\left(V^\ell(y^u+a_n+h^\ell(\tin))-h_0(\tin)\right)
\Big|
\end{eqnarray*}
by (\ref{glivenko-cantelli}) and the definitions of $[s^\ell,s^u]$ and $[y^\ell,y^u]$. Further, we can bound the last sum by
\begin{eqnarray*}
 &&{} \frac 1n\sum_{i=1}^n\Big| F_0\left(V^u(y^u+a_n+h^u(\tin))-h_0(\tin)\right)\\
&&{}\qquad-F_0\left(\Lambda_0(\Lambda_\t^{-1}(y^u+a_n+h^u(\tin)))-h_0(\tin)\right)\Big|\\
&&{}+ \frac 1n\sum_{i=1}^n\Big| F_0\left(\Lambda_0(\Lambda_\t^{-1}(y^u+a_n+h^\ell(\tin)))-h_0(\tin)\right)\\
&&\qquad{}-F_0\left(V^\ell(y^u+a_n+h^\ell(\tin))-h_0(\tin)\right)\Big|\\
&&{}+ \frac 1n\sum_{i=1}^n\Big| F_0\left(\Lambda_0(\Lambda_\t^{-1}(y^u+a_n+h^u(\tin)))-h_0(\tin)\right)\\
&&\qquad{}
-F_0\left(\Lambda_0(\Lambda_\t^{-1}(y^u+a_n+h^\ell(\tin)))-h_0(\tin)\right)\Big|\\
 &\leq& 3\epsilon
\end{eqnarray*}
using the construction of brackets above (note that $\|V^u-\Lambda_0\circ\Lambda_\t^{-1}\|_\infty\leq \delta$, $\|\Lambda_0\circ\Lambda_\t^{-1}-V^\ell\|_\infty\leq \delta$, $\|h^u-h^\ell\|_\infty\leq\gamma$   and recall assumptions \ref{F2} and \ref{L2}).

Thus $\sup_{\t\in\Theta}\sup_{s\in [0,1]\atop y\in C}|H_{n,\t}^{(1)}(y,s)|$ can be bounded by $O(\epsilon)+o(1)$ plus a finite maximum over the absolute value of terms 
\begin{eqnarray*}
\frac 1n\sum_{i=1}^n \Big(I\{\Lambda_0(Y_{i,n})\leq V^u(y^u+a_n+h^u({\tin}))\} -  \mathbb{E}[I\{\Lambda_0(Y_{i,n})\leq V^u(y^u+a_n+h^u({\tin}))\}]\Big)
\end{eqnarray*}
and
\begin{eqnarray*}
\frac 1n\sum_{i=1}^n \Big(I\{\Lambda_0(Y_{i,n})\leq V^\ell(y^\ell+a_n+h^\ell({\tin}))\} -  \mathbb{E}[I\{\Lambda_0(Y_{i,n})\leq V^\ell(y^\ell+a_n+h^\ell({\tin}))\}]\Big).
\end{eqnarray*}
However, those converge to zero in probability by a simple application of Chebychev's inequality. 

This completes the proof of $\sup_{\t\in\Theta}\|H_{n,\t}^{(1)}\|=o_P(1)$ and thus of the lemma.  \hfill $\Box$

\begin{lemma}\label{lem1} Under the assumptions of Theorem \ref{theo} (ii),
$$\sup_{\t\in\Theta} \|G(\t,h_{\t})-G(\t,\hat h_{\t})\| =o_P(1).$$
\end{lemma}

{\bf Proof.}
The proof is analogous to the proof of Lemma \ref{lem1-random}.
 \hfill $\Box$

\begin{lemma}\label{lem3} Under the assumptions of Theorem \ref{theo} (ii),
$$\sup_{\t\in\Theta}\|\tilde G_n(\t,\hat h_{\t})-G(\t,\hat h_{ \t})\| =o_P(1).$$
\end{lemma}

{\bf Proof.}
According to assumption \ref{N1}  it suffices to show
$$\sup_{\t\in\Theta}\sup_{s\in [0,1]\atop y\in C} |\tilde G_n(\t,\hat h_{\t})(y,s)-G(\t,\hat h_{ \t})(y,s)| =o_P(1).$$
Recalling the definitions of $\tilde G_n$ in (\ref{G_n-tilde}) and $G$ in (\ref{G-random}) we only consider the first sum and first integral, respectively. 
It holds by the mean value theorem for integration
\beq
&&\Bigg|\frac 1n\sum_{i=1}^n F_0\left(\Lambda_0(\Lt^{-1}(y+\hat h_\t(\tin)))-h_0(\tin)\right)I\{x_{i,n}\leq s\}\\
&&{}-\int F_0\left(\Lambda_0(\Lt^{-1}(y+\hat h_\t(x)))\!-\!h_0(x)\right)I\{x\leq s\} f_X(x)\,dx\Bigg|\\
&=&\bigg|\sum_{i=1}^n \Big(\frac 1nF_0\left(\Lambda_0(\Lt^{-1}(y+\hat h_\t(\tin)))-h_0(\tin)\right)\\
&&{}-\int_{x_{i-1,n}}^{x_{i,n}}F_0\left(\Lambda_0(\Lt^{-1}(y+\hat h_\t(x)))-h_0(x)\right)f_X(x)dx\Big)I\{x_{i,n}\leq s\}\\
&&{}-\int_{\max\{x_{i,n}| x_{i,n}\leq s\}}^sF_0\left(\Lambda_0(\Lt^{-1}(y+\hat h_\t(x)))-h_0(x)\right)f_X(x)\,dx\bigg|\\
&\leq& \sum_{i=1}^n \bigg| \frac 1nF_0\left(\Lambda_0(\Lt^{-1}(y+\hat h_\t(\tin)))-h_0(\tin)\right)\\
&&{}-F_0\left(\Lambda_0(\Lt^{-1}(y+\hat h_\t(\xi_{i,n})))-h_0(\xi_{i,n})\right) f_X(\xi_{i,n})(x_{i,n}-x_{i-1,n})\bigg|\\
&&{}+O(\bar\Delta_n)
\eeq
for some $\xi_{i,n}\in[x_{i-1,n},x_{i,n}]$. Now the assertion follows from assumption \ref{X1''}, uniform continuity of $F_0$ and of $\Lambda_0\circ\Lambda_\t^{-1}$ (uniformly in $\t$)
and from 
$$|\hat h_\t(\tin)-\hat h_\t(\xi_{i,n})|\leq \|\hat h_\t-h_\t\|_\infty+|\Lt(\Lambda_0^{-1}( h_0(x_{i,n}))-\Lt(\Lambda_0^{-1}( h_0(\xi_{i,n})))|$$ 
in connection with Lemma \ref{lemma-smoothing} and assumptions \ref{H1}, \ref{L1}. 
 \hfill $\Box$

%

\section{Identifiability of the model in the fixed design case}\label{ident-proof-fixed}

To prove identifiability in the case of deterministic covariates as in Remark \ref{identifiability-fixed} one starts similarly to the proof in section \ref{ident-proof} of the appendix (main paper) with the cdf of $\varepsilon_{i,n}(\vartheta_1)=\Lambda_{\vartheta_1}(Y_{i,n})-h_{\vartheta_1}(x_{i,n})$ in $y$ to obtain that 
$H^{-1}(y+H(h_{\vartheta_0}(x_{i,n})))- h_{\vartheta_0}(x_{i,n})$
does not depend on $x_{i,n}$ for $y\in (-\infty,0]$. Due to continuity of the functions and $\bar{\Delta}_n\to 0$ one obtains that  $H^{-1}(y+H(h_{\vartheta_0}(x)))- h_{\vartheta_0}(x)$
does not depend on $x\in[0,1]$ for $y\in (-\infty,0]$. The remainder of the proof is as in section \ref{ident-proof}.

\section*{References}

\begin{description}

\item Drees, H., Neumeyer, N. and Selk, L. (2018). Estimation and hypotheses tests in boundary regression models. \textit{Bernoulli}, to appear.  \\
http://www.bernoulli-society.org/index.php/publications/bernoulli-journal

\end{description}

\section{Figures and Tables}

\setcounter{figure}{4}
\setcounter{table}{5}
 
\begin{figure}[h]
\center\epsfig{file=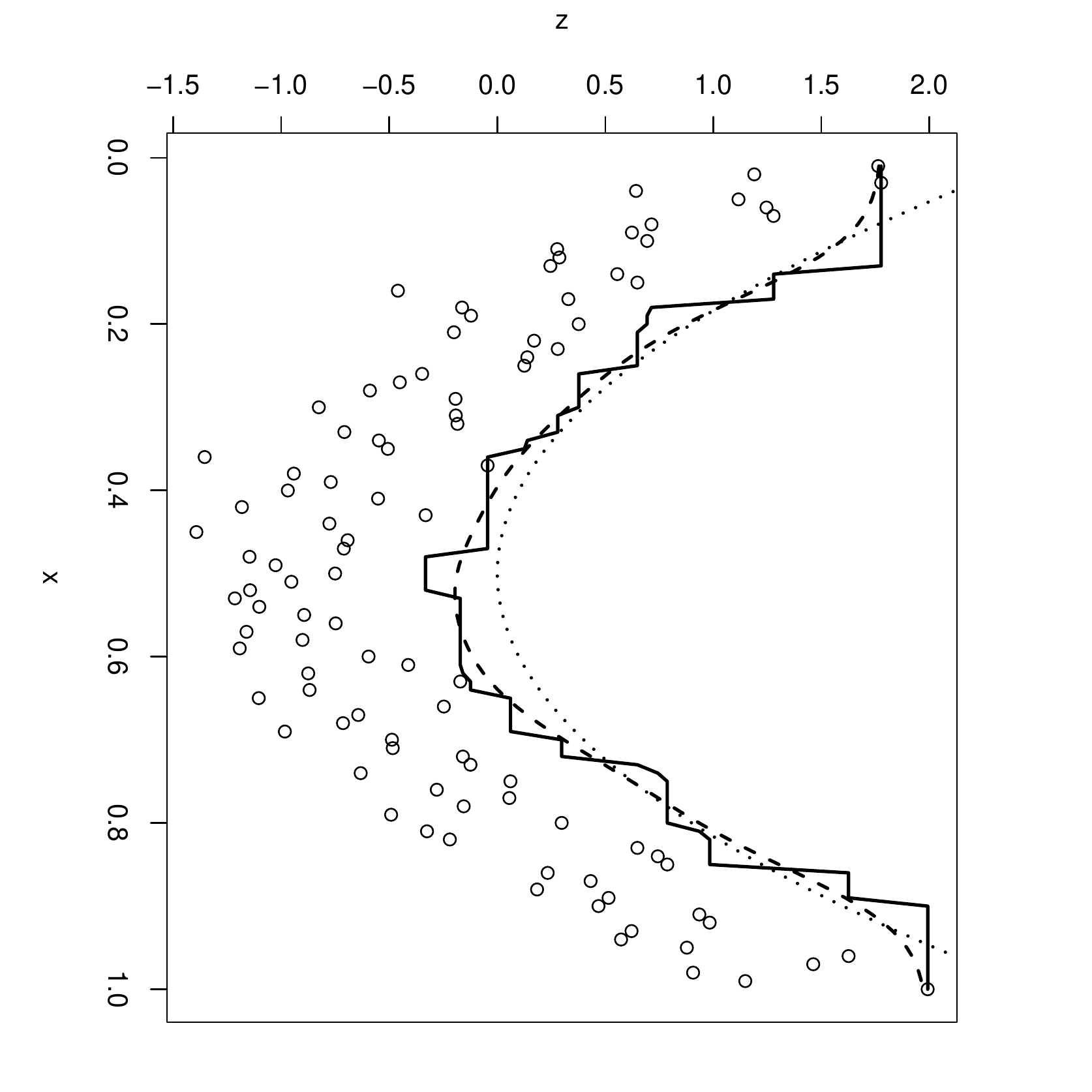,angle=90,width=11cm}
\caption{Data corresponding to the model in Figure \ref{graphic1}. The true curve is dotted, while the local constant estimator is given by the solid line and the smoothed estimator (with bandwidths $b_n=n^{-1/3}$ and $a_n=b_n/2$) by the dashed line. 
}\label{graphic3}
\end{figure}

\begin{figure}[h]
\center\epsfig{file=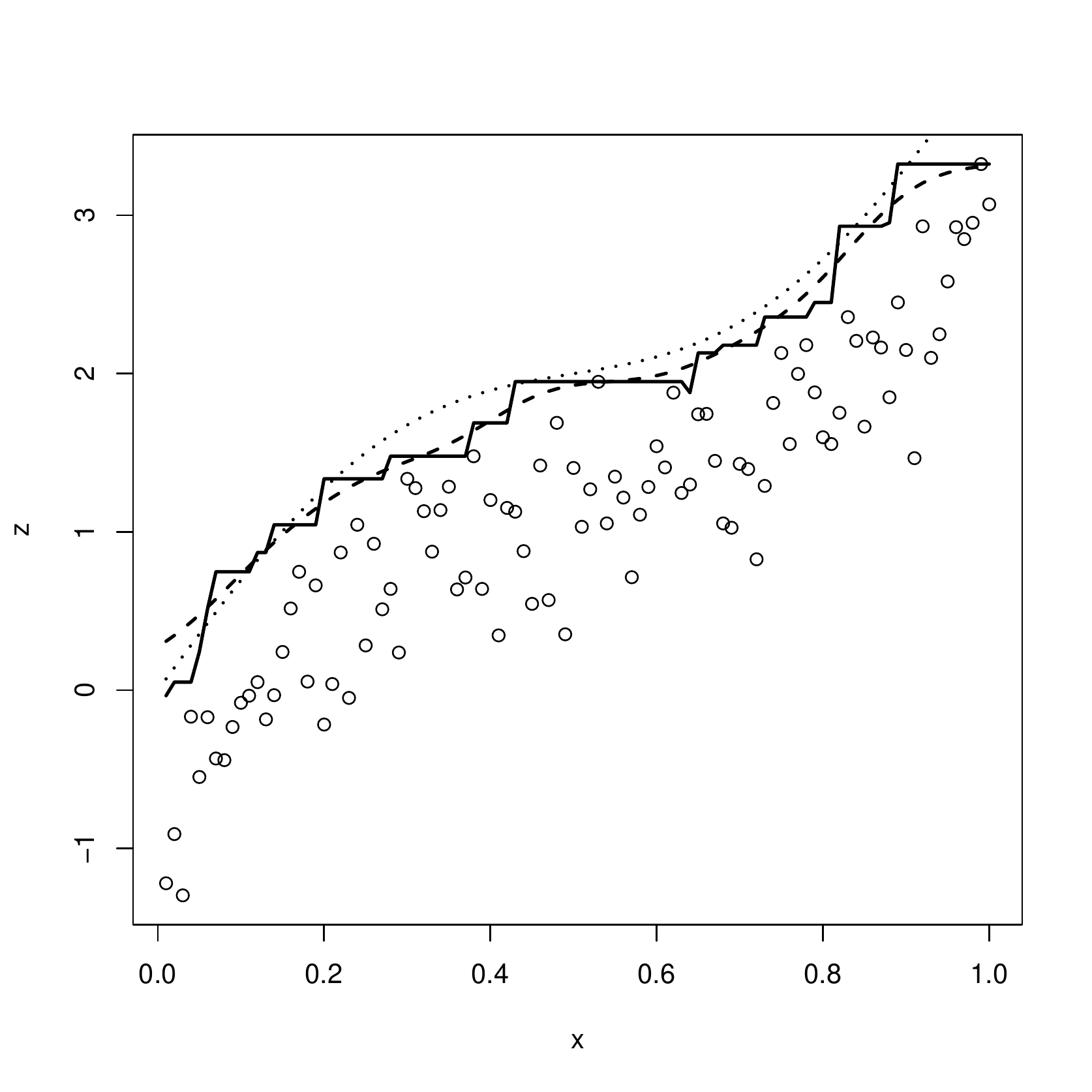,width=11cm}
\caption{The setting is similar to Figure \ref{graphic3} with model from Figure \ref{graphic2} in the main paper. 
}\label{graphic4}
\end{figure}

\begin{table}
\begin{tabular}{|l|c|c|c|c|}
			\hline
			$n=50$& TKS & TCM & TKSCM & TCMKS \\
			\hline
			$\t_0=0$
			& 0.192 0.197 (0.102) & \textbf{0.001} 0.009 (\underline{0.037}) & -0.038 -0.139 (0.118) & 0.225 0.208 (0.118)
			\\
			\hline
			$\t_0=0.5$
			& 0.778 0.691 (0.191) & \textbf{0.378} 0.402 (\underline{0.092}) & 0.239 0.302 (0.410) & 0.858 0.798 (0.274)
			\\
			\hline
			$\t_0=1$
			& 1.290 1.340 (0.233) & \textbf{0.728} 0.741 (\underline{0.232}) & 0.388 0.264 (1.000) & 1.370 1.350 (0.308)
			\\
			\hline
			$\t_0=1.5$
			& \textbf{1.750} 1.790 (\underline{0.195}) & 1.160 1.290 (0.368) & 0.507 0.292 (1.810) & 1.790 1.790 (0.222)
			\\
			\hline
			$\t_0=2$
			& 1.940 2.060 (0.201) & 1.590 1.750 (0.478) & 0.585 0.424 (2.880) & \textbf{1.970} 2.060 (\underline{0.141})
			\\
			\hline
		\end{tabular}
		
		\bigskip
		
		\begin{tabular}{|l|c|c|c|c|}
			\hline
			$n=100$ & TKS & TCM & TKSCM & TCMKS \\
			\hline
			$\t_0=0$
			& \textbf{0.017} 0.018 (0.037) & 0.080 0.079 (\underline{0.014}) & 0.061 0.073 (0.022) & -0.020 -0.004 (0.020)
			\\
			\hline
			$\t_0=0.5$
			& \textbf{0.496} 0.517 (\underline{0.028}) & 0.338 0.346 (0.042) & 0.516 0.578 (0.080) & 0.521 0.548 (0.032) 
			\\
			\hline
			$\t_0=1$
			& 0.973 0.979 (\underline{0.044}) & 0.745 0.745 (0.092) & 0.906 1.050 (0.225) & \textbf{1.030} 1.020 (0.054)
			\\
			\hline
			$\t_0=1.5$
			& 1.480 1.460 (\underline{0.059}) & 1.210 1.230 (0.123) & 1.310 1.510 (0.412) & \textbf{1.510} 1.490 (0.060)
			\\
			\hline
			$\t_0=2$
			& \textbf{1.960} 2.000 (0.059) & 1.690 1.740 (0.144) & 1.550 1.860 (0.822) & 1.920 1.940 (\underline{0.058})
			\\
			\hline
		\end{tabular}
		\captionof{table}  {
		\it Mean, median and MISE for Model (\ref{simu2}) for $n=50$ and $n=100$ with $a_n=b_n/2$. }\label{Table3}
\end{table}

\begin{table}	
		\begin{tabular}{|l|c|c|c|c|}
			\hline
			$n=50$& TKS & TCM & TKSCM & TCMKS \\
			\hline
			$\t_0=0$
			& 0.156 0.167 (0.076) & -0.050 -0.062 (\underline{0.031}) & \textbf{0.022} 0.103 (0.095) & 0.191 0.198 (0.0,86)
			\\
			\hline
			$\t_0=0.5$
			& 0.713 0.646 (0.130) & \textbf{0.324} 0.336 (\underline{0.088}) & 0.268 0.407 (0.348) & 0.781 0.695 (0.197)
			\\
			\hline
			$\t_0=1$
			& 1.260 1.310 (\underline{0.191}) & 0.655 0.646 (0.242) & 0.447 0.511 (0.919) & \textbf{1.330} 1.350 (0.258)
			\\
			\hline
			$\t_0=1.5$
			& \textbf{1.720} 1.780 (0.188) & 1.100 1.180 (0.365) & 0.619 0.559 (1.660) & \textbf{1.720} 1.780 (\underline{0.177})
			\\
			\hline
			$\t_0=2$
			& \textbf{1.970} 2.060 (0.141) & 1.550 1.660 (0.442) & 0.726 0.619 (2.630) & 1.960 2.060 (\underline{0.111})
			\\
			\hline
		\end{tabular}
		
		\bigskip
		
		\begin{tabular}{|l|c|c|c|c|}
			\hline
			$n=100$ & TKS & TCM & TKSCM & TCMKS \\
			\hline
			$\t_0=0$
			& \textbf{0.001} 0.050 (0.044) & 0.129 0.128 (0.023) & 0.028 0.037 (\underline{0.016}) & -0.014 -0.042 (\underline{0.015})
			\\
			\hline
			$\t_0=0.5$
			& 0.467 0.474 (\underline{0.033}) & 0.282 0.287 (0.063) & \textbf{0.497} 0.533 (0.057) & 0.481 0.486 (\underline{0.034}) 
			\\
			\hline
			$\t_0=1$
			& 0.934 0.942 (\underline{0.043}) & 0.674 0.649 (0.130) & 0.878 0.999 (0.190) & \textbf{0.965} 0.960 (0.049)
			\\
			\hline
			$\t_0=1.5$
			& 1.420 1.390 (0.056) & 1.120 1.130 (0.185) & 1.320 1.500 (0.336) & \textbf{1.440} 1.400 (\underline{0.053})
			\\
			\hline
			$\t_0=2$
			& \textbf{1.910} 1.920 (\underline{0.071}) & 1.590 1.610 (0.228) & 1.560 1.850 (0.790) & 1.850 1.790 (0.079)
			\\
			\hline
		\end{tabular}
		\captionof{table}  {
		\it Mean, median and MISE for Model (\ref{simu2}) for $n=50$ and $n=100$ with $a_n=b_n/20$. }\label{Table4}
\end{table}		

\begin{table}
\begin{tabular}{|l|c|c|c|c|}
			\hline
			$n=50$& TKS & TCM & TKSCM & TCMKS \\
			\hline
			$\t_0=0$
			& 0.141 0.010 (0.309) & \textbf{0.020} 0.012 (\underline{0.008}) & -0.056 -0.053 (0.016) & 0.137 0.028 (0.260)
			\\
			\hline
			$\t_0=0.5$
			& 0.519 0.549 (0.039) & \textbf{0.518} 0.506 (\underline{0.023}) & 0.521 0.532 (0.040) & 0.546 0.574 (0.038)
			\\
			\hline
			$\t_0=1$
			& 1.010 1.010 (0.085) & \textbf{1.000} 0.996 (\underline{0.040}) & 0.996 0.998 (0.071) & 1.040 1.030 (0.077)
			\\
			\hline
			$\t_0=1.5$
			& 1.530 1.530 (0.125) & \textbf{1.500} 1.490 (\underline{0.066}) & \textbf{1.500} 1.510 (0.113) & 1.550 1.570 (0.110)
			\\
			\hline
			$\t_0=2$
			& 1.960 2.060 (0.118) & \textbf{2.010} 2.040 (\underline{0.069}) & 1.950 2.000 (0.156) & 1.970 2.050 (0.093)
			\\
			\hline
		\end{tabular}
		
		\bigskip
		
		\begin{tabular}{|l|c|c|c|c|}
			\hline
			$n=100$ & TKS & TCM & TKSCM & TCMKS \\
			\hline
			$\t_0=0$
			& 0.019 0.009 (0.022) & \textbf{0.006} 0.000 (\underline{0.004}) & 0.043 0.038 (0.007) & -0.014 -0.007 (0.013)
			\\
			\hline
			$\t_0=0.5$
			& 0.522 0.524 (0.023) & \textbf{0.505} 0.498 (\underline{0.013}) & 0.562 0.555 (0.020) & 0.528 0.524 (0.022)
			\\
			\hline
			$\t_0=1$
			& 1.030 1.030 (0.042) & \textbf{1.010} 1.000 (\underline{0.021}) & 1.080 1.080 (0.038) & 1.030 1.020 (0.042)
			\\
			\hline
			$\t_0=1.5$
			& 1.550 1.550 (0.061) & \textbf{1.510} 1.510 (\underline{0.030}) & 1.600 1.590 (0.055) & 1.550 1.550 (0.061)
			\\
			\hline
			$\t_0=2$
			& 2.040 2.060 (0.066) & \textbf{2.000} 2.000 (\underline{0.037}) & 2.070 2.070 (0.061) & 2.020 2.050 (0.058)
			\\
			\hline
		\end{tabular}
	\captionof{table}  {
		\it Mean, median and MISE for Model (\ref{simu3}) for $n=50$ and $n=100$ with $a_n=b_n/2$. }\label{Table5}
\end{table}

\begin{table}
			\begin{tabular}{|l|c|c|c|c|}
			\hline
			$n=50$& TKS & TCM & TKSCM & TCMKS \\
			\hline
			$\t_0=0$
			& 0.097 0.031 (0.223) & \textbf{0.005} 0.016 (\underline{0.008}) & -0.037 -0.030 (0.013) & -0.087 -0.001 (0.174)
			\\
			\hline
			$\t_0=0.5$
			& 0.487 0.506 (0.039) & 0.479 0.462 (\underline{0.021}) & \textbf{0.506} 0.508 (0.036) & 0.514 0.522 (0.035)
			\\
			\hline
			$\t_0=1$
			& 0.976 0.978 (0.092) & 0.965 0.962 (\underline{0.044}) & \textbf{0.984} 0.997 (0.074) & 1.020 1.020 (0.078)
			\\
			\hline
			$\t_0=1.5$
			& \textbf{1.499} 1.469 (0.120) & 1.440 1.430 (\underline{0.063}) & 1.450 1.470 (0.119) & 1.530 1.500 (0105)
			\\
			\hline
			$\t_0=2$
			& 1.920 1.990 (0.105) & \textbf{1.960} 1.960 (\underline{0.069}) & 1.940 1.940 (0.127) & 1.930 1.970 (0.086)
			\\
			\hline
		\end{tabular}
		
		\bigskip
		
		\begin{tabular}{|l|c|c|c|c|}
			\hline
			$n=100$ & TKS & TCM & TKSCM & TCMKS \\
			\hline
			$\t_0=0$
			& 0.017 0.004 (0.016) & \textbf{0.004} 0.000 (\underline{0.004}) & 0.042 0.039 (0.007) & -0.010 -0.007 (0.007)
			\\
			\hline
			$\t_0=0.5$
			& 0.530 0.537 (0.021) & \textbf{0.507} 0.497 (\underline{0.011}) & 0.563 0.554 (0.019) & 0.534 0.539 (0.021)
			\\
			\hline
			$\t_0=1$
			& 1.020 1.020 (0.042) & \textbf{1.000} 1.000 (\underline{0.020}) & 1.080 1.070 (0.035) & 1.020 1.010 (0.039)
			\\
			\hline
			$\t_0=1.5$
			& 1.550 1.560 (0.064) & \textbf{1.510} 1.510 (\underline{0.031}) & 1.600 1.600 (0.054) & 1.560 1.550 (0.064)
			\\
			\hline
			$\t_0=2$
			& 2.050 2.060 (0.069) & \textbf{2.020} 2.040 (\underline{0.041}) & 2.090 2.100 (0.064) & 2.030 2.060 (0.059)
			\\
			\hline
		\end{tabular}
		\captionof{table}  {
		\it Mean, median and MISE for Model (\ref{simu3}) for $n=50$ and $n=100$ with $a_n=b_n/20$. }\label{Table6}
\end{table}

\begin{table}
	\begin{center}
				\begin{tabular}{|l|c|c|c|}
			\hline
			Method & Pearson & Kendall & Spearman \\
			\hline
			Original data
			& -0.634 & -0.456 & -0.612   
			\\
			\hline
			True parameter $\vartheta_0$
			& 0.001 & 0.001 & 0.001
			\\
			\hline
			TKS
			& 0.001 & 0.001 & 0.001
			\\
			\hline
			TCM
			& 0.009 & 0.006 & 0.008
			\\
			\hline
			TKSCM
			& 0.005 & 0.002 & 0.004
			\\
			\hline
			TCMKS & 0.002 & 0.001 & 0.001
			\\
			\hline
		\end{tabular}
		\captionof{table}  {
		\it Pearson's, Kendall's and Spearman's correlation coefficients (the average over $1000$ iterations) between the covariates and the errors for the model \eqref{simu3} when $n=100$. The first line corresponds to the correlations for the original data while the second line is for the true transformation parameter ($\vartheta_0=0.5$). The last four lines correspond to the correlations for each estimator.}\label{TableCor1paper}
\end{center}
\end{table}

\end{appendix}

\end{document}